\documentclass[11pt]{amsart}

\usepackage{amsmath}
\usepackage{amssymb}
\usepackage[all]{xy}
\usepackage{sty1}
\usepackage{thsty}

\makeatletter
\@addtoreset{equation}{section}
\makeatother

\title[Measurable Schur Multipliers]{Measurable Schur Multipliers
and Completely Bounded Multipliers of the Fourier Algebras}

\author{Nico Spronk}

\date{}

\def\tr{\mathrm{tr}}
\def\diag{\mathrm{diag}}

\begin{document}
\maketitle

\newtheorem{slicemaps}{Proposition}[section]
\newtheorem{ehtensprodchar}[slicemaps]{Theorem}
\newtheorem{ehtensprodchar1}[slicemaps]{Corollary}
\newtheorem{fhiseh}[slicemaps]{Theorem}

\newtheorem{classicschurmult}{Theorem}[section]
\newtheorem{classicschurmult1}[classicschurmult]{Corollary}
\newtheorem{schurmult}[classicschurmult]{Theorem}
\newtheorem{cschurmult}[classicschurmult]{Proposition}
\newtheorem{exthaagerup2}[classicschurmult]{Proposition}
\newtheorem{exthaagerup3}[classicschurmult]{Proposition}

\newtheorem{fsalgccba}{Proposition}[section]
\newtheorem{mpliersismults}[fsalgccba]{Proposition}
\newtheorem{cbmclosure}[fsalgccba]{Proposition}
\newtheorem{fmultschar}[fsalgccba]{Proposition}
\newtheorem{fsalgchar}[fsalgccba]{Proposition}
\newtheorem{fsalgchar1}[fsalgccba]{Corollary}

\newtheorem{invschurmult1}{Lemma}[section]
\newtheorem{invschurmult2}[invschurmult1]{Theorem}
\newtheorem{haagerup}[invschurmult1]{Theorem}
\newtheorem{haagerup0}[invschurmult1]{Corollary}
\newtheorem{haagerup1}[invschurmult1]{Corollary}
\newtheorem{cbmgiswap1}[invschurmult1]{Corollary}

\newtheorem{jsiginv}{Proposition}[section]
\newtheorem{jsiginv1}[jsiginv]{Theorem}
\newtheorem{jsiginv2}[jsiginv]{Corollary}
\newtheorem{jsiginv3}[jsiginv]{Corollary}
\newtheorem{cbmpredual}[jsiginv]{Theorem}
\newtheorem{cbmpredual1}[jsiginv]{Corollary}
\newtheorem{vmapinjective}[jsiginv]{Proposition}
\newtheorem{vmapinjective1}[jsiginv]{Corollary}
\newtheorem{qmapsurj}[jsiginv]{Proposition}
\newtheorem{gamenchar}[jsiginv]{Theorem}
\newtheorem{simprob}[jsiginv]{Theorem}

\footnote{{\it Date}: \today.


This work was supported at various stages by an NSERC PGS B and
by Ontario Graduate Scholarships.
}

\begin{abstract}
Let $G$ be a locally compact group, $\mathrm{L}^p(G)$ be
the usual L$^p$-space for $1\leq p\leq\infty$
and $\falg$ be the Fourier algebra of $G$.  Our goal is to
study, in a new abstract context, 
the completely bounded multipliers of $\falg$, which we denote
$\cbmg$.  We show that $\cbmg$ can be characterised as the ``invariant part''
of the space of (completely) bounded normal $\blinftyg$-bimodule maps
on $\bdop{\bltwog}$, the space of bounded operator on $\bltwog$.
In doing this we develop a function theoretic description of the
normal $\blinfty{X,\mu}$-bimodule maps on $\bdop{\bltwo{X,\mu}}$, which we 
denote by $\schur{X,\mu}$, and name the {\it measurable Schur multipliers} 
of $(X,\mu)$.  Our approach leads to many new results, some of which generalise
results hitherto known only for certain classes of groups.  Those results which
we develop here are a uniform approach to obtaining the functorial properties
of $\cbmg$, and a concrete description of a standard predual
of $\cbmg$.  
\end{abstract}

\section{Introduction}

Let $G$ be a locally compact group, and $\falg$ the Fourier algebra
of $G$ from \cite{eymard}.

In \cite{herz74}, Herz defined the space of {\it Herz-Schur multipliers},
$\mathrm{B}_2(G)$.  He claimed this space to coincide with the space
he denoted $\fV(G)$ in his article \cite{herz}, but offered no proof of 
this fact.  The {\it completely bounded multipliers} of the Fourier 
algebra, $\cbmg$, were defined by De Canneire and Haagerup in 
\cite{decanniereh}.  In was shown by Bo\.{z}eko and Fendler,
\cite{bozejkof}, that $\cbmg=\mathrm{B}_2(G)$, isometrically.
Unfortunately, this result relies on unpublished work of Gilbert 
\cite{gilbert}.
However, there is a proof of the fact that $\cbmg=\mathrm{B}_2(G)$
due to Jolissaint \cite{jolissaint}, which relies on the representation
theorem for completely bounded maps applied to the reduced C*-algebra,
$\cstarrg$.

In this paper we prove an analogue of the fact that $\cbmg=\mathrm{B}_2(G)$.
In fact, our result obtains that fact in the case that $G$ is a 
discrete group, and improves upon it in the respect that we obtain a
complete isometry, where both $\cbmg$ and $\mathrm{B}_2(G)$ are given natural 
operator space structures.  Moreover, we show that the natural
module action of $\cbmg$ on the von Neumann algebra generated by
the left regular representation, $\vng$, extends completely isometrically
to a normal action of $\cbmg$ on all of $\bdop{\bltwog}$.  In fact, 
the image of $\cbmg$ in the (completely) bounded normal maps
on $\bdop{\bltwog}$ consists of
$\blinftyg$-bimodule maps, and can be characterised
as a certain intrinsic ``invariant part'' amongst those maps.
We feel that our approach is interesting in and of itself since it makes
systematic use of modern operator space theoretic techniques, pioneered
by Blecher, Effros, Haagerup, Paulsen, Ruan and Smith, amongst others.
It thus gives a starting point
for obtaining representations of completely multipliers of Kac algebras,
after \cite{krausr} or of multipliers of locally compact quantum groups after
Vaes et al (see \cite{vaes}).

To support our approach, in Section \ref{sec:schurmults}
we describe the space of 
{\it measurable Schur multipliers},
$\schur{X,\mu}$, for any suitable measure space $(X,\mu)$.  We develop
the theory of $\schur{X,\mu}$ fully and systematically.  Our main tools
for doing this are the results of Smith \cite{smith91} and
Blecher and Smith \cite{blechers92} which characterise the
normal completely bounded maps on the space of all bounded operators on 
a Hilbert space, which are module maps over a von Neumann algebra.
We make great use of the weak* Haagerup tensor product
of two von Neumann algebras, which is
described in \cite{blechers92}.  
In particular, we note that $\schur{X,\mu}$ obtains a natural
operator space structure.
Subordinate to the theory of measurable Schur multipliers, we develop
a theory of ``continuous'' Schur multipliers, which is based on the
theory of the extended Haagerup tensor product of Effros and Ruan
in \cite{effrosrP} (see also \cite{effroskr}).

Letting $m$ now denote the left Haar measure on our locally compact group $G$,
we let $\schur{G}=\schur{G,m}$.  We define the ``invariant part'',
$\schurinv{G}$, of $\schur{G}$.  Our space $\schurinv{G}$ is essentially
the space $\fV(G)$ of \cite{herz}.  Our main result is that
$\schurinv{G}\cong\cbmg$ completely isometrically.  Thus we obtain
a natural extension of the module action of $\cbmg$ on
$\vng$ to $\bdop{\bltwog}$.  This is the content of Section \ref{sec:invschur}.
(It should be noted that the isometric identification $\schurinv{G}\cong\cbmg$
has been previously known to Haagerup \cite{haagerup80}.  This fact 
was discovered by the author after the present article was completed.
However, it is not clear how the methods of that article can be adapted to
give a complete isometry.)

As an immediate application of our results, we have a systematic approach
to determining the functorial properties of $\cbmg$.  In addition,
we can identify some circumstances in which a subspace of the Fourier-Stieltjes
algebra, $\fsalg$, acts completely isometrically as multipliers of $\falg$.
Another application we obtain is a concrete description of the predual
of $\cbmg$, $\queg$, whose existence is recognised in \cite{decanniereh}.
In the case where $G$ is discrete, Pisier \cite{pisier98}
implicitly uses the structure of $\queg$ in his partial solution
to the similarity problem for representations of groups.
These applications are studied in Section \ref{sec:applications}.

We note that there are further applications of our results and methods
which are not explored in this note.  As an example, the author and
L.\ Turowska, in \cite{spronkt}, use some of the methods of the current note
to show that if $G$ is compact, then $\falg$ naturally imbeds into
$\cont{G}\tens^\gam\cont{G}$ (projective tensor product of the continuous
functions on $G$ with itself).  This is used to generalise, to arbitrary
compact groups, results on
parallel spectral synthesis due to Varopoulos \cite{varopoulos} for 
compact Abelian groups; and to generalise some results on operator synthesis
due to Froelich \cite{froelich}.

We should also point out that if $G$ is Abelian, there are related
results to ours by St{\o}rmer \cite{stormer}.  In that article it is shown that
there is an isometric representation of the measure algebra, $\measg$,
in $\wbdop{\bltwog}$, which extends convolution on $\blinftyg$.  
This is generalised to arbitrary locally compact 
groups by Ghahramani \cite{ghahramani}.  Neufang \cite{neufang} has shown 
that this isometric representation has its range in $\wcbop{\bltwog}$.
Some of the connections between these results and those of the present
article are being explored by Neufang and the author \cite{neufangs}.

Sections \ref{sec:haageruptp} and \ref{sec:cbmults} are included to support
the exposition of this article, but include results of independant interest.
In Section \ref{sec:haageruptp}, a uniform approach is developed for defining
the weak* and extended Haagerup tensor products.  In Section \ref{sec:cbmults}
some basic results about completely bounded multipliers are reviewed, and
a theorem of Walter \cite{walter89} is improved.

\smallskip
This article represents a significant portion of the author's doctoral thesis
\cite{spronkTh}.  The author would like to express his gratitude to
his advisor, Brian Forrest, for support and guidance throughout his
tenure as a graduate student at the University of Waterloo.  His gratitude
extends to other members of the faculty there, including
Ken Davidson, Kathryn Hare and Andu Nica.  The author would also like
to acknowledge his fellow graduate student
Peter Wood, and the excellent staff in the Pure Math department.

\subsection{Notation}
\label{ssec:notation}

If $\fX$ is a complex normed vector space, we will let $\ball{\fX}$
denote the closed unit ball, $\fX^*$ denote the space of continuous 
linear functionals on $\fX$, i.e.\ the dual space, and 
$\bdop{\fX}$ denote the normed algebra of bounded linear operators on $\fX$.
If $T\in\bdop{\fX}$, we let $T^*:\fX^*\to\fX^*$ denote the adjoint operator.
If $\hil$ is a Hilbert space and $x\in\bdop{\hil}$, 
we let $x^*:\hil\to\hil$ denote its Hilbertian adjoint.

If $\fX$ is a dual space, with established predual $\fX_*$, we let
$\wbdop{\fX}$ denote the space of linear operators on $\fX$
which are continuous in the $\sig(\fX,\fX_*)$ topology.  Hence if
$\fM$ is a von Neumann algebra, $\wbdop{\fM}$ denotes the algebra of
normal operators on $\fM$.

If $X$ is a locally compact Hausdorff space, let $\contb{X}$ denote
the C*-algebra of continuous complex-valued bounded functions on $X$,
with norm $\unorm{\vphi}=\sup\{|\vphi(x)|:x\in X\}$.  We will let
$\conto{X}$ denote the norm closure of the space of compactly supported
continuous complex-valued functions on $X$.

Our standard reference for operator spaces is \cite{effrosrB}.  However
we were also influenced by the points of view in \cite{blecher} and 
\cite{blecherp}.  There is a good summary of operator spaces in
\cite{runde}.  

Given an operator space $\fV$, we let $\matn{\fV}$ denote the space of
$n\cross n$ matrices with entries in $\fV$.
We also let $\cbop{\fV}$ denote the algebra of
completely bounded linear maps on $\fV$, i.e.\ maps $T:\fV\to\fV$
for which the amplifications, $\nlift{T}:\matn{\fV}\to\matn{\fV}$,
satisfy $\cbnorm{T}=\sup\{\norm{\nlift{T}}:n\in\En\}<+\infty$.
We will assign $\cbop{\fV}$ and $\fV^*$ the operator space
from \cite{blecher} or \cite[Sec.\ 3.2]{effrosrB}.

If $\fX$ is a normed vector space, we let $\max\fX$ denote the
maximal operator space with underlying normed space $\fX$, and
we denote by $\min\fX$ the minimal one.

We will make substantial use of tensor products here.  If $\fV$ and
$\fW$ are complete operator spaces, let $\fV\ptens\fW$ denote their {\it 
operator projective tensor product}, and $\fV\htens\fW$ their {\it Haagerup
tensor product}.  We will often write $\fV\pptens\fW$, respectively
$\fV\phtens\fW$, for the algebraic tensor product $\fV\otimes\fW$, but with
the operator space structure given by the operator projective tensor norms,
respectively Haagerup, tensor norms.  For Banach spaces $\fX$ and $\fY$,
$\fX\tens^\gam\fY$ will denote their (Banach space) {\it projective tensor 
product}.  If $\fM$ and $\fN$ are von Neumann algebras
we will denote the {\it von Neumann tensor
product} by $\fM\vntens\fN$.

\section{The Weak* and Extended Haagerup Tensor Products}
\label{sec:haageruptp}

In this article we will make extensive use of the weak* Haagerup
tensor product of \cite{blechers92} and the extended Haagerup tensor
product of \cite{effrosrP} (see also \cite{effroskr}).  However, for our 
applications we will not want to assume that we can work with operator spaces
over a separable Hilbert space.  Moreover, we want a representation of
the extended Haagerup tensor product modeled after one of the weak*
Haagerup tensor product in \cite[Theo.\ 3.1]{blechers92}.  Hence we will
briskly develop a theory of these tensor products in a form that fits
our needs.

First we need a description of infinite matrices over an operator space.
Let $\fV$ be an operator space and $I_0$ and $J_0$ be finite index sets.
The set of $I_0\cross J_0$-matrices, $\mat{I_0,J_0}{\fV}$, can be normed
in an unambiguous way with no regard to having an ordering on $I_0$ or $J_0$;
see \cite[(2.1.5)]{effrosrB}.  If $I$ and $J$ are infinite index sets,
let $M_{I,J}(\fV)$ denote the set of $I\cross J$ matrices with entries in
$\fV$.  If $I_0$ and $J_0$ are finite subsets of $I$ and $J$, and 
$v\in M_{I,J}(\fV)$, let $v_{I_0,J_0}$ be the matrix in $\mat{I_0,J_0}{\fV}$
given by restricting to the indices from $I_0$ and $J_0$.  We let
\[
\mat{I,J}{\fV}=\cbrac{v\in M_{I,J}(\fV):\norm{v}=\sup_{\scriptsize
\begin{matrix} I_0\subset I, J_0\subset J \\ \text{finite}\end{matrix}
}\norm{v_{I_0,J_0}}<+\infty}.
\]
As in \cite[Sec.\ 10.1]{effrosrB}, this is an operator space.  For example,
if $\hil$ is a Hilbert space, then $\mat{I,J}{\bdop{\hil}}\cong
\bdop{\hil^J,\hil^I}$, in the usual way.  In particular,
$\smat{I,J}=\mat{I,J}{\Cee}\cong\bdop{\ltwo{J},\ltwo{I}}$.

If $\fW$ is another operator space and $S\in\cbop{\fV,\fW}$, then the 
amplification $\lift{I,J}{S}:\mat{I,J}{\fV}\to M_{I,J}(\fW)$ can be easily
seen to have its range in $\mat{I,J}{\fW}$ and to be completely bounded
with $\norm{\lift{I,J}{S}}=\cbnorm{\lift{I,J}{S}}=\cbnorm{S}$.
Moreover, $\lift{I,J}{S}$ is a complete isometry if $S$ is.  However,
if $S$ is a complete quotient map, $\lift{I,J}{S}$ may not be a
(complete) quotient map;  it is such, though, if
$\fV$ and $\fW$ are operator dual spaces
and $S$ is the adjoint of a complete isometry between their preduals.

Usually, a matrix $v\iin\mat{I,J}{\fV}$ will be written
$v=\sbrac{v_{ij}}$.  However, for row matrices $v\iin\mat{1,J}{\fV}$
and column matrices $w\iin\mat{I,1}{\fV}$ we will often write
$v=\sbrac{v_j}$ and $w=\sbrac{w_i}$, i.e.\ $v_j=v_{1j}$ and $w_i=w_{i1}$
for each $i$ and $j$.  

\subsection{The Fubini-Haagerup Tensor Product}

The Fubini-Haagerup tensor product, which we define below,
is developed in the same way as the weak* Haagerup tensor product
is in \cite{blechers92}, and is general enough to give us a usable
form of the extended Haagerup tensor product of \cite{effrosrP}.

Let us begin by recalling the theory of the weak* Haagerup tensor product,
but in the context we require.  Let $\hil$ be a Hilbert space of Hilbertian
dimension $\dim\hil$ and $I$ be an index set of cardinality
$|I|=(\dim\hil)^2$.  Squaring accommodates the case that $\dim\hil$ is finite.
Following arguments from \cite[Theo.\ 3.1]{smith91},
which may also be found in  \cite{haagerup80},
we have that if $S\in\wcbop{\bdop{\hil}}$ then there exist
$\sbrac{v_i}\iin\mat{1,I}{\bdop{\hil}}$ and 
$\sbrac{w_i}\iin\mat{I,1}{\bdop{\hil}}$ such that for any $a$ in
$\bdop{\hil}$
\[
Sa=\sum_{i\in I} v_iaw_i=
\begin{bmatrix} v_{i_1} & v_{i_2} & \cdots \end{bmatrix}a
\begin{bmatrix} w_{i_1} \\ w_{i_2} \\ \vdots \end{bmatrix}
\]
where we well-order $I=\{i_1,i_2,\dots\}$ to represent matrix
multiplication, and, moreover
\[
\cbnorm{S}=\norm{\sbrac{v_i}}\norm{\sbrac{w_i}}=
\norm{\sum_{i\in I} v_iv_i^*}^{1/2}\norm{\sum_{i\in I} w_i^*w_i}^{1/2}.
\]
The size of the index set $I$ may be deduced by tracing through the proofs
the representation theorem for completely bounded maps found in
\cite[Theo.\ 2.7]{paulsen84} and of
Stinespring's Theorem \cite{stinespring}; or see
\cite[Sec.\ 5.2]{effrosrB} for these results together.
The aforementioned arguments can be adapted to see that if
$\sbrac{S_{ij}}\in\matn{\wcbop{\bdop{\hil}}}\cong
\wcbop{\bdop{\hil},\matn{\bdop{\hil}}}$ then there are matrices
$\sbrac{v_{i\iota}}\iin\mat{n,I}{\bdop{\hil}}$ and
$\sbrac{w_{\iota j}}\iin\mat{I,n}{\bdop{\hil}}$ such that
\[
S_{ij}a=\sum_{\iota\in I} v_{i\iota}aw_{\iota j}\quad
\aand\quad \cbnorm{S}=\norm{\sbrac{v_{i\iota}}}\norm{\sbrac{w_{\iota j}}}.
\]
Then, following \cite{blechers92}, we define
\begin{equation}\label{eq:whtensprod}
\bdop{\hil}\whtens\bdop{\hil}=
\cbrac{\sum_{i\in I} v_i\tens w_i:
\begin{matrix} \sbrac{v_i}\in\mat{1,I}{\bdop{\hil}}\aand \\
\sbrac{w_i}\in\mat{I,1}{\bdop{\hil}}\end{matrix}}.
\end{equation}
Here, two sums of elementary tensors,  $\sum_{i\in I} v_i\tens w_i$
and $\sum_{i\in I} v_i'\tens w_i'$ will be regarded to be the same 
if $\sum_{i\in I} v_iaw_i=\sum_{i\in I} v_i'aw_i'$ for every $a$ in
$\bdop{\hil}$.  
The norm on $\matn{\bdop{\hil}\whtens\bdop{\hil}}$ is given by
\begin{equation}
\label{eq:matrixwhnorm}
\whnorm{u}=\min\cbrac{\norm{v}\norm{w}:\begin{matrix} 
v\in\mat{n,I}{\bdop{\hil}},\; w\in\mat{I,n}{\bdop{\hil}} \\ 
\aand u=v\odot w\end{matrix}}
\end{equation}
where for $v=\sbrac{v_{i\iota}}\iin\mat{n,I}{\bdop{\hil}}$ and
$w=\sbrac{w_{\iota j}}\iin\mat{I,n}{\bdop{\hil}}$, $v\odot w$ has
$ij$th component given by
\[
(v\odot w)_{ij}=\sum_{\iota\in I} v_{i\iota}\tens w_{\iota j}.
\]
If $u=\sum_{i\in I} v_i\tens w_i$ in 
$\bdop{\hil}\whtens\bdop{\hil}$, define $T_u:\bdop{\hil}\to\bdop{\hil}$
by
\begin{equation}\label{eq:tmap}
T_ua=\sum_{i\in I} v_ia w_i.
\end{equation}
Then, by our construction of $\bdop{\hil}\whtens\bdop{\hil}$, 
$u\mapsto T_u$ is a surjective complete isometry from
$\bdop{\hil}\whtens\bdop{\hil}$ to $\wcbop{\bdop{\hil}}$.
As in \cite{blechers92}, we have that
$\bdop{\hil}\whtens\bdop{\hil}$ is a dual operator space.
We also note that it can be deduced from 
\cite{smith91} that $\bdop{\hil}\htens\bdop{\hil}$ imbeds completely
isometrically into $\bdop{\hil}\whtens\bdop{\hil}$.

The {\it left} and {\it right slice maps}, $L_F$ and $R_F$ below,
will be our major tool for identifying useful subspaces of
$\bdop{\hil}\whtens\bdop{\hil}$.

\begin{slicemaps}\label{prop:slicemaps}
If $F\in\bdop{\hil}^*$ and $u=\sum_{i\in I} v_i\tens w_i\iin
\bdop{\hil}\whtens\bdop{\hil}$, the sums
\begin{equation}\label{eq:slicemap}
L_Fu=\sum_{i\in I} F(v_i)w_i\quad\aand\quad
R_Fu=\sum_{i\in I} F(w_i)v_i.
\end{equation}
converge in norm.  Moreover, $L_F$ and $R_F$ define completely bounded
maps from $\bdop{\hil}\whtens\bdop{\hil}$ to $\bdop{\hil}$.
\end{slicemaps} 

\proof We will show this only for $L_F$.

To see that $\sum_{i\in I} F(v_i)w_i$ converges in norm, we simply
note that since $F$ is automatically completely bounded we have
that $\lift{1,I}{F}\sbrac{v_i}=\sbrac{F(v_i)}\in\smat{1,I}$,
where $\smat{1,I}$ is the row Hilbert space $\ltwo{I}_r$, so the
coefficients $F(v_i)$ are square summable.

To see that $L_F$ defines a a completely bounded linear map, we must
first check that it is well-defined.  If $\xi,\eta$ and $\zeta$ are vectors
in $\hil$, let $\ome_{\xi,\eta}$ denote the usual vector functional
and $\zeta\tens\eta^*$ the usual rank 1 operator.  If
$\sbrac{v_i}\iin\mat{1,I}{\bdop{\hil}}$ and 
$\sbrac{w_i}\iin\mat{I,1}{\bdop{\hil}}$ are such that
$\sum_{i\in I}v_iaw_i=0$ for every $a\iin\bdop{\hil}$, then,
letting $u=\sum_{i\in I} v_i\tens w_i$ we have that
\[
\sbrac{R_{\ome_{\xi,\eta}}u}\zeta 
= \sum_{i\in I} \inprod{w_i\xi}{\eta}v_i\zeta 
= \sum_{i\in I}\inprod{\xi}{w_i^*\eta}v_i\zeta 
= \sum_{i\in I}v_i(\zeta\tens\eta^*)w_i\xi
=0\mult\xi=0
\]
from which it follows that $R_{\ome_{\xi,\eta}}u=0$.  Then, we have that
$\ome_{\xi,\eta}(L_Fu)=F(R_{\ome_{\xi,\eta}}u)=0$ for each $\xi$ and
$\eta$ in $\hil$, so $L_Fu=0$. 

The complete boundedness of $L_F$ is now easy to check.  Indeed,
if $u=v\odot w$ in $\matn{\bdop{\hil}\whtens\bdop{\hil}}$, as in
(\ref{eq:matrixwhnorm}), with $\whnorm{u}=\norm{v}\norm{w}$, then
\[
\norm{\nlift{L_F}u}=\norm{\brac{\lift{n,I}{F}v}w}\leq
\norm{F}\norm{v}\norm{w}=\norm{F}\whnorm{u}
\]
so $\cbnorm{L_F}\leq\norm{F}$. 

(Thanks go to D.\ Blecher for pointing out to the author that it is
not {\it a priori} clear that $L_F\aand R_F$ are well-defined.
The proof here uses the idea of \cite[Prop.\ 3.7 (i)]{blechers92}.)
\endpf

\smallskip
Now let $\fV$ and $\fW$ be closed subspaces of $\bdop{\hil}$.
We define the {\it Fubini-Haagerup tensor product} by
\[
\fV\fhtens\fW=\cbrac{u\in\bdop{\hil}\whtens\bdop{\hil}:
\begin{matrix} L_Fu\in\fW\aand R_Fu\in\fV \\
\text{for each }F\iin\bdop{\hil}^*\end{matrix}}.
\]
It is clear that $\fV\fhtens\fW$ is a norm-closed subspace
of $\bdop{\hil}\whtens\bdop{\hil}$.
It follows from the theorem below that $\fV\fhtens\fW$ can be defined
for abstract operator spaces, independently of their completely
isometric representations in $\bdop{\hil}$.

\begin{ehtensprodchar}\label{theo:ehtensprodchar}
Let $\fV$ and $\fW$ be closed subspaces of $\bdop{\hil}$ and
$u\in\bdop{\hil}\whtens\bdop{\hil}$.  Then the following are equivalent:

{\bf (i)} $u\in\fV\fhtens\fW$.

{\bf (ii)} $R_\ome u\in\fV$ and $L_\ome u\in\fW$ for each
$\ome\iin\bdop{\hil}_*$.

{\bf (iii)} There exist $v\iin\mat{1,I}{\fV}$ and $w\iin\mat{I,1}{\fW}$
such that
\[
u=v\odot w\quad\aand\quad\whnorm{u}=\norm{v}\norm{w}.
\]
In particular, we have that if $u\in\matn{\fV\fhtens\fW}$, then
\begin{equation}\label{eq:ehelement}
\whnorm{u}=\min\cbrac{\norm{v}\norm{w}:\begin{matrix} v\in\mat{n,I}{\fV},
w\in\mat{I,n}{\fW} \\ \aand u=v\odot w\end{matrix}}.
\end{equation}
\end{ehtensprodchar}

\proof That (i) implies (ii) is trivial.
If (ii) holds, then we have that $R_\ome u\in\fV$ and $L_\ome u\in\fW$ for 
each vector functional $\ome=\ome_{\xi,\eta}$.  Then
the proof follows identically that of \cite[Theo.\ 4.1]{smith91},
adapted to accommodate uncountable index sets;
we observe that the definition of {\it strong independence} 
\cite[Def.\ 3.2]{smith91} is suitable for uncountable
column matrices of operators.
Finally, if (iii) holds then (i) follows from Proposition 
\ref{prop:slicemaps}.
\endpf


\begin{ehtensprodchar1}
\label{cor:ehtensprodchar1}
If $\fV,\fW,\fX$ and $\fY$ are all operator spaces, $S\in\cbop{\fV,\fX}$
and $T\in\cbop{\fW,\fY}$, then the map $S\tens T:\fV\tens\fW\to\fX\tens\fY$
extends to a completely bounded map, again denoted $S\tens T$, from
$\fV\fhtens\fW$ to $\fX\fhtens\fY$.  If $S$ and $T$ are each complete 
isometries, then so too is $S\tens T$.
\end{ehtensprodchar1}

\proof It should first be checked that $S\tens T$ is well-defined on
$\fV\fhtens\fW$.  If $u=0$ in  $\fV\fhtens\fW$, then we must have
that $(S\tens T)u=0$ since for each $F\iin\fX^*$ we have that
$\sbrac{L_F\comp(S\tens T)}u=TL_{F{\scriptstyle\circ} S}u=T0=0$, 
where the fact that
$L_{F\comp S}u=0$ follows from Proposition \ref{prop:slicemaps}
and the Arveson-Wittstock Extension Theorem.  Similarly,
$\sbrac{R_{F'}\comp(S\tens T)}u=0$ for each $F'\iin\fY^*$. 
It then follows easily from (\ref{eq:matrixwhnorm})
that $\cbnorm{S\tens T}\leq\cbnorm{S}\cbnorm{T}$.

That $S\tens T$ is a complete isometry when $S$ and $T$ each are, also follows
from (\ref{eq:matrixwhnorm}).  \endpf

\smallskip
It is also immediate from Theorem \ref{theo:ehtensprodchar} that
if $\fV$ and $\fW$ are dual operator spaces then $\fV\fhtens\fW=
\fV\whtens\fW$, the {\it weak* Haagerup tensor product};
see \cite[Theo.\ 3.1]{blechers92}.
Hereafter we prefer to use the notation $\whtens$ in this case.
As in \cite[Theo.\ 3.2]{blechers92}, 
if $\fV$ and $\fW$ have operator preduals $\fV_*$ and $\fW_*$, 
then $\fV\whtens\fW\cong\brac{\fV_*\htens\fW_*}^*$ via the dual pairing
\begin{equation}\label{eq:whdualpair}
\dpair{\sum_{i\in I}v_i\tens w_i}{\nu\tens\omega}
=\sum_{i\in I}\nu(v_i)\omega(w_i).
\end{equation}

If $\fV$ and $\fW$ are any pair of operator spaces, the {\it extended
Haagerup tensor product} is defined in \cite{effrosrP} by
\[
\fV\ehtens\fW=\cbrac{B\in\brac{\fV^*\htens\fW^*}^*:
\begin{matrix} f\mapsto B(f\tens g_0)\aand g\mapsto B(f_0\tens g) \\
\text{are weak* continuous for any} \\
g_0\iin\fW^*\aand f_0\iin\fV^* \end{matrix}}
\]

\begin{fhiseh}
For any pair of operator spaces $\fV$ and $\fW$,
$\fV\ehtens\fW\cong\fV\fhtens\fW$ completely isometrically.
\end{fhiseh}

\proof Both $\fV\ehtens\fW$ and $\fV\fhtens\fW$ naturally inject
completely isometrically into
$\brac{\fV^*\htens\fW^*}^*\cong\fV^{**}\whtens\fW^{**}$,
with the same range. \endpf

\smallskip
We will now dispense with the notation $\fV\fhtens\fW$, preferring to use
$\fV\ehtens\fW$ in its place, unless both of $\fV$ and $\fW$ are dual spaces.

\section{Schur Multipliers}
\label{sec:schurmults}

\subsection{Classical Schur Multipliers}

We motivate our study in Subsection \ref{ssec:measschurmult} by
exposing the theory of classical Schur multipliers in our context.

Let us first fix an index set $I$ and consider the Hilbert space $\ltwo{I}$.
It will be convenient to identify $\bdop{\ltwo{I}}$ with
$\smat{I,I}$, the space of $I\cross I$ scalar matrices representing
bounded operators on $\ltwo{I}$, with respect to the canonical
ortho-normal basis $\{\del_i\}_{i\in I}$.

A function $u:I\cross I\to\Cee$ is called a {\it Schur multiplier}
if $\sbrac{u(i,j)\alp_{ij}}\in\smat{I,I}$ for every
$\sbrac{\alp_{ij}}\iin\smat{I,I}$.  An application of the Closed Graph
Theorem shows that the operator
\begin{equation}\label{eq:schurmults}
S_u:\smat{I,I}\to\smat{I,I}\quad\text{ given by }\quad
S_u\sbrac{\alp_{ij}}=\sbrac{u(i,j)\alp_{ij}}
\end{equation}
is bounded.  We denote the space of Schur multipliers by $\schur{I}$
and let $\scnorm{u}=\norm{S_u}$ for each $u\iin\schur{I}$.

The von Neumann algebra $\linfty{I}$ of bounded functions on $I$ can
be naturally identified with $\{M_\vphi:\vphi\in\linfty{I}\}$, the
algebra of multiplication operators on $\ltwo{I}$:  if
$\vphi\in\linfty{I}$, then
$M_\vphi\del_i=\vphi(i)\del_i$.  We thus identify
$\linfty{I}$ with the algebra of diagonal matrices in $\smat{I,I}$.

In the notation of the previous section, we have that
\[
\linfty{I}\whtens\linfty{I}=
\cbrac{\sum_{k\in K}\vphi_k\tens\psi_k:
\unorm{\sum_{k\in K}|\vphi_k|^2}^{1/2}\unorm{\sum_{k\in K}|\psi_k|^2}^{1/2}
<+\infty}.
\]
Here, $K$ is an index set chosen so that $|K|=|I|^2$, where
the squaring accommodates the case that $I$ is finite.
We know from Theorem \ref{theo:ehtensprodchar} that this space has norm
given by
\[
\whnorm{v}=\min\cbrac{
\unorm{\sum_{k\in K}|\vphi_k|^2}^{1/2}\unorm{\sum_{k\in K}|\psi_k|^2}^{1/2}:
v=\sum_{k\in K}\vphi_k\tens\psi_k}.
\]
We have, then, that the operator $T_v:\smat{I,I}\to\smat{I,I}$ takes the form
\begin{equation}\label{eq:diagmult}
T_v\sbrac{\alp_{ij}}=\sum_{k\in K} 
\diag(\vphi_k)\sbrac{\alp_{ij}}\diag(\psi_k)
\end{equation}
where $\diag(\vphi)$ is the diagonal matrix with diagonal entries 
$\diag(\vphi)_{ii}=\vphi(i)$.

The following result is actually \cite[Prop.\ 1.1]{pisier95}, but with
a different proof.

\begin{classicschurmult}\label{theo:classicschurmult}
A function $u:I\cross I\to\Cee$ is in $\schur{I}$ if and only if
there is a $v=\sum_{k\in K}\vphi_k\tens\psi_k
\iin\linfty{I}\whtens\linfty{I}$ such that $u=v$, i.e.\
\[
u(i,j)=\sum_{k\in K}\vphi_k(i)\psi_k(j)
\]
for each $(i,j)\iin I\cross I$.  Moreover, $S_u=T_v$ in this case,
and hence $\scnorm{u}=\whnorm{v}$.
\end{classicschurmult}

\proof First, if $v=\sum_{k\in K}\vphi_k\tens\psi_k
\iin\linfty{I}\whtens\linfty{I}$, write
\[
v(i,j)=\sum_{k\in K}\vphi_k(i)\psi_k(j)
\]
for each $(i,j)\iin I\cross I$.  
Then if $\sbrac{\alp_{ij}}\in\smatfin{I,I}$, i.e.\ $\sbrac{\alp_{ij}}$ has
only finitely many non-zero entries, it follows from
(\ref{eq:diagmult}) that
\begin{equation}\label{eq:diagmultcalc}
T_v\sbrac{\alp_{ij}}=\sum_{k\in K}\sbrac{\vphi_k(i)\alp_{ij}\psi_k(j)}
=\sbrac{\sum_{k\in K}\vphi_k(i)\alp_{ij}\psi_k(j)}
=\sbrac{v(i,j)\alp_{ij}}.
\end{equation}
Since $T_v$ is weak*-weak* continuous on $\smat{I,I}$ and 
$\smatfin{I,I}$ is weak*-dense in $\smat{I,I}$, 
(\ref{eq:diagmultcalc}) holds for
any $\sbrac{\alp_{ij}}\iin\smat{I,I}$.  In particular, we see that
$v$, qua function on $I\cross I$, is in $\schur{I}$ and 
with $T_v=S_v$.

Now if $u\in\schur{I}$, let $\check{u}(i,j)=u(j,i)$.  Then
$\scnorm{\check{u}}=\scnorm{u}$.  Indeed the transpose
map $x\mapsto\transpose{x}$ is an isometry on $\smat{I,I}$ and
$S_{\check{u}}=\transpose{\brac{S_u\transpose{x}}}$ for each $x$.
Let $\stmat{I,I}$ denote the $I\cross I$ trace class matrices,
so $\stmat{I,I}^*\cong\smat{I,I}$ via $\dpair{t}{x}=\tr(tx)$.
We note that $S_{\check{u}}(\stmatf{I,I})\subset\stmatf{I,I}$
and that $S_u=\brac{S_{\check{u}}|_{\stmatf{I,I}}}^*$.  
This implies that $S_{\check{u}}$ is bounded on $\stmatf{I,I}$
and hence extends to $\stmat{I,I}$, with 
$S_u=\brac{S_{\check{u}}|_{\stmat{I,I}}}^*$.  Hence
$S_u$ is weak*-weak* continuous on $\smat{I,I}$.  Moreover,
$S_u$ is clearly an $\linfty{I}$-bimodule map.
Since $\linfty{I}$ is locally cyclic on $\ltwo{I}$ (and, in fact,
cyclic if $I$ is countable), $S_u$ is completely bounded
by \cite[Sec.\ 2]{smith91}.  Hence we find that
$S_u=T_v$ for some $v$ in $\linfty{I}\whtens\linfty{I}$,
by \cite[Theo.\ 3.1]{smith91}.  
It then follows from above that
$v=u$ and that $\whnorm{v}=\norm{T_v}=\norm{S_u}=\scnorm{u}$. \endpf

\smallskip
We can use Theorem \ref{theo:classicschurmult} to deduce the following
well-known result of Grothendieck.  See \cite[Prop.\ 1.1]{pisier95}
or \cite[Theo.\ 5.1]{pisierB}, for example.

\begin{classicschurmult1}\label{cor:classicschurmult1}
A function $u:I\cross I\to\Cee$ is in $\schur{I}$ if and only if
there is a Hilbert space $\hil$, and bounded functions
$\xi,\eta:I\to\hil$ such that
\[
u(i,j)=\inprod{\xi(i)}{\eta(j)}
\]
for each $(i,j)\iin I\cross I$.  We write $u=u_{\xi,\eta}$, in this case.
Moreover, $\scnorm{u}=\inf\{\unorm{\xi}\unorm{\eta}:u=u_{\xi,\eta}\}$.
\end{classicschurmult1}

\proof  If $u\in\schur{I}$, then $u\in\linfty{I}\whtens\linfty{I}$
and we may write $u=\sum_{k\in K}\vphi_k\tens\psi_k$ where
$\scnorm{u}=\unorm{\sum_{k\in K}|\vphi_k|^2}^{1/2}
\unorm{\sum_{k\in K}|\psi_k|^2}^{1/2}$.  Let $\xi,\eta:I\to\ltwo{K}$
be given by
\begin{equation}\label{eq:vectorfromfunc}
\xi(i)=\brac{\vphi_k(i)}_{k\in K}\quad\aand\quad
\eta(i)=\brac{\wbar{\psi_k(i)}}_{k\in K}.
\end{equation}
Then $\unorm{\xi}=\unorm{\sum_{k\in K}|\vphi_k|^2}^{1/2}$ and
$\unorm{\eta}=\unorm{\sum_{k\in K}|\psi_k|^2}^{1/2}$, and
$u=u_{\xi,\eta}$.

Conversely, let $u=u_{\xi,\eta}$ where $\xi,\eta:I\to\hil$ are bounded 
functions.  Let $\{e_k\}_{k\in K}$ be an ortho-normal basis for $\hil$.
Then for each $k$ in $K$, let $\vphi_k,\psi_k:I\to\Cee$ be given by
\[
\vphi_k(i)=\inprod{\xi(i)}{e_k}\quad\aand\quad
\psi_k(i)=\inprod{e_k}{\eta(i)}.
\]
Then it follows from Parseval's formula that 
$\sum_{k\in K}|\vphi_k(i)|^2=\norm{\xi(i)}^2$, so
$\unorm{\sum_{k\in K}|\vphi_k|^2}^{1/2}=\unorm{\xi}$.
Similarly, $\unorm{\sum_{k\in K}|\psi_k|^2}^{1/2}=\unorm{\eta}$ and
we get that $u=\sum_{k\in K}\vphi_k\tens\psi_k$.  \endpf

\subsection{Measurable Schur Multipliers}
\label{ssec:measschurmult}

For this subsection we let $(X,\mu)$ denote a measure space,
with positive measure $\mu$,
which satisfies the Radon-Nikodym Theorem so that $\blone{X,\mu}^*\cong
\blinfty{X,\mu}$.  This is satisfied whenever $\mu$ is $\sig$-finite or
if $\mu$ is a regular measure on the Borel sets of a locally compact space
$X$ (see \cite[III.12]{hewittrI} for the latter).  We note that since 
$(X,\mu)$ is
not necessarily $\sig$-finite, elements of $\blinfty{X,\mu}$ can be
identified as functions only up to locally null sets.  A subset $N$ of $X$
is called {\it $\mu$-locally null} if $\mu(N\cap E)=0$ for any $\mu$-measurable
subset $E$ of $X$ such that $\mu(E)<+\infty$.  Thus the norm on 
$\blinfty{X,\mu}$ is given by
\begin{align*}
\unorm{\vphi}&=\esssup_{x\in X}|\vphi(x)| \\
&=\inf\big\{c\in\Ree^+:\{ x\in X:|\vphi(x)|>c\}\text{ is locally null}\big\}.
\end{align*}

We will want spaces of Hilbert space valued $\mathrm{L}^\infty$-functions.
For our purposes it is sufficient to consider only weakly measurable
such functions in the sense of \cite{diestelu}.  We outline the
construction of such spaces to clarify our notation.
If $\hil$ is a Hilbert space, let $\blinftyhil{X,\mu}{\hil}_0$ denote
the set of functions $\xi:X\to\hil$ such that

{\bf (i)} for all $\eta\iin\hil$, $x\mapsto\inprod{\xi(x)}{\eta}$ is
           measurable, and

{\bf (ii)} $\unorm{\xi}=\esssup_{x\in X}\norm{\xi(x)}<+\infty$.

\noindent Note that if $\hil$ is separable, then (ii) follows from

{\bf (ii')} $\sup\{\unorm{\inprod{\xi(\cdot)}{\eta}}:
            \eta\in\ball{\hil}\}<+\infty$

\noindent by a straightforward application of the 
Uniform Boundedness Principle.

It is evident that $\blinftyhil{X,\mu}{\hil}_0$ is a linear space with
semi-norm $\unorm{\cdot}$.  Hence if $\blinftyhil{X,\mu}{\hil}$ denotes
the quotient space $\blinftyhil{X,\mu}{\hil}_0/\{\xi:\unorm{\xi}=0\}$,
it is a normed space, which can be verified to be complete.  As is
the case with $\blinfty{X,\mu}$, we will consider elements of
$\blinftyhil{X,\mu}{\hil}$ as functions, defined $\mu$-locally almost
everywhere.

A subset $E$ of $X\cross X$ is called {\it $\mu\cross\mu$-locally marginally 
null}, if $E\subset(N_1\cross X)\cup(X\cross N_2)$ for some locally null sets
$N_1\aand N_2$ in $X$.  We will be interested in describing spaces of functions
defined only up to locally marginally null sets.  We note that if $X$ is
$\sig$-finite for $\mu$, then every $\mu\cross\mu$-locally marginally 
null is a $\mu\cross\mu$-null set.  However, the converse does not hold:
consider $X=[0,1]$, the unit interval in $\Ree$, 
with $\mu$ being Lebesgue measure and $E=\{(s,s):s\in[0,1]\}$.
For convenience,
we will abbreviate ``(locally) marginally almost every'' to (l.)m.a.e.\ in the 
sequel.

If $\xi,\eta\in\blinftyhil{X,\mu}{\hil}$, define $u_{\xi,\eta}$ for
$\mu\cross\mu$-l.m.a.e.\ $(x,y)\iin X\cross X$ by
\begin{equation}
\label{eq:schurfunc}
u_{\xi,\eta}(x,y)=\inprod{\xi(x)}{\eta(y)}.
\end{equation}
Note that if $\xi=\xi'\aand\eta=\eta'$, except on $\mu$-locally null
sets $N_\xi$ and $N_\eta$, then $u_{\xi,\eta}=u_{\xi',\eta'}$ except
on $(N_\xi\cross X)\cup(X\cross N_\eta)$,
so there is no ambiguity in defining $u_{\xi,\eta}$.  Let
\begin{equation}\label{eq:mschur}
\schur{X,\mu}=\{u_{\xi,\eta}:\xi,\eta\in\blinftyhil{X,\mu}{\hil}
\text{ for some Hilbert space }\hil\}.
\end{equation}
By using Hilbert space direct sums and tensor products, it is easy to
see that $\schur{X,\mu}$ is an algebra.  In fact, it can be further
verified that it is a Banach algebra under the norm
\begin{equation}\label{eq:scnorm}
\scnorm{u}=\inf\{\unorm{\xi}\unorm{\eta}:u=u_{\xi,\eta}\}.
\end{equation}
for $u\iin\schur{X,\mu}$.  However, since this fact follows from 
Theorem \ref{theo:schurmult} {\it infra}, we will not prove it directly.

For $p=1,2\oor\infty$, let $\mathrm{L}^p=\blp{X,\mu}$ below.

The space $\tee{X,\mu}=\abltwo\tens^\gam\abltwo$ is the predual of 
$\bdop{\abltwo}$ via the identification
\begin{equation}
\label{eq:tdualisb}
\dpair{S}{f\tens g}=\inprod{Sf}{\wbar{g}}=\int_X Sf(x)g(x)d\mu(x).
\end{equation}
Elements of $\tee{X,\mu}$ may be regarded as functions, defined up to
marginally null sets:  if $\ome=\sum_{k=1}^\infty f_k\tens g_k\in\tee{X,\mu}$,
then let
\begin{equation}\label{eq:teefunc}
\ome(x,y)=\sum_{k=1}^\infty f_k(y) g_k(x)
\end{equation}
for m.a.e.\ $(x,y)\iin X\cross X$.  See \cite[2.2.7]{arveson}.  Note that the 
order of $x$ and $y$ above is purposeful, and leads to technical 
simplifications in the sequel.  A function $u:X\cross X\to\Cee$ is said to be a
{\it multiplier} of $\tee{X,\mu}$ if
\begin{equation}
\label{eq:multstee}
m_u\ome(x,y)=u(x,y)f(y)g(x)
\end{equation}
defines an element of $\tee{X,\mu}$ for each elementary tensor $\ome=f\tens g$ 
in $\abltwo\tens_\gam\abltwo$ and does so in such a way that
$\norm{m_u f\tens g}_\gam\leq C\norm{f}_2\norm{g}_2$ for a fixed constant 
$C>0$.  Hence $m_u$ extends to a bounded linear map on $\tee{X,\mu}$.
It is easy to see that two multipliers $u$ and $u'$ satisfy
$m_u=m_{u'}$ if and only if $u=u'$ locally marginally almost everywhere.
The linear space of ($\mu\cross\mu$-l.m.a.e.\ equivalence classes of)
such functions will be denoted $\mults\tee{X,\mu}$.

Let 
\begin{equation}
\label{eq:hilfac}
\hilfac{\ablone}{\ablinfty}=\cbrac{T:\ablone\to\ablinfty:
\begin{matrix} T=T_2\comp T_1\wwhere T_1\in\bdop{\ablone,\hil} \\
\aand  T_2\in\bdop{\hil,\ablinfty}, \\
\text{for some Hilbert space }\hil \end{matrix}
}.
\end{equation}
The norm on $\hilfac{\ablone}{\ablinfty}$ is given by
\begin{equation*}
\hfnorm{T}=\inf\{\norm{T_2}\norm{T_1}:T=T_2\comp T_1 
\text{ as in (\ref{eq:hilfac})}\}.
\end{equation*}
We note that since $\ablone=\max\ablone$  and $\ablinfty=\min\ablinfty$,
$\hilfac{\ablone}{\ablinfty}=\colfac{\ablone}{\ablinfty}$ with the same
norm, where
$\colfac{\ablone}{\ablinfty}$ is the space of operators which factor
as completely bounded maps through a column space.  See \cite[13.3]{effrosrB}
or \cite{effrosr91}.  Hence, by \cite[Theo.\ 5.3]{effrosr91}
$\hilfac{\ablone}{\ablinfty}\cong(\ablone\htens\ablone)^*$ via
\begin{equation}\label{eq:blonehd}
\dpair{T}{f\tens g}=\dpair{Tf}{g}=\int_X Tf(x)g(x)d\mu(x)
\end{equation}
and thus $\hilfac{\ablone}{\ablinfty}$ attains the dual
operator space structure.  We note that the matrix norms on 
$\hilfac{\ablone}{\ablinfty}$
may be described independently of this duality (see \cite[13.3]{effrosrB}).

Since $(X,\mu)$ is assumed to satisfy the Radon-Nikodym Theorem, we have
that $(\ablone)^*\cong\ablinfty$.  Moreover, we have a natural identification
of $\ablinfty$ with $\{M_\vphi:\vphi\in\ablinfty\}$, the maximal Abelian
subalgebra of $\bdop{\abltwo}$ consisting of multiplication operators,
where $M_\vphi f=\vphi f$ for $f\iin\abltwo$.  Thus, following
(\ref{eq:whtensprod}), we obtain the space
\[
\ablinfty\whtens\ablinfty=\cbrac{\sum_{i\in I}\vphi_i\tens\psi_i:
\unorm{\sum_{i\in I}|\vphi_i|^2}\unorm{\sum_{i\in I}|\psi_i|^2}<+\infty}
\]
with norm
\[
\scnorm{u}=\min\cbrac{
\unorm{\sum_{i\in I}|\vphi_i|^2}^{1/2}\unorm{\sum_{i\in I}|\psi_i|^2}^{1/2}:
u=\sum_{i\in I}\vphi_i\tens\psi_i}.
\]
It will be convenient to represent the duality $(\ablone\htens\ablone)^*\cong
\ablinfty\whtens\ablinfty$ given in (\ref{eq:whdualpair}):
\begin{equation}
\label{eq:lonehdual}
\dpair{\sum_{i\in I}\vphi_i\tens\psi_i}{f\tens g}
=\sum_{i\in I}\brac{\int_X \vphi_i f d\mu}\brac{ \int_X \psi_i g d\mu}
\end{equation}
for $f,g\iin\ablone$.  
By \cite[Cor.\ 3.8]{blechers92}, $\ablinfty\whtens\ablinfty$
injects naturally into $\ablinfty\vntens\ablinfty
\cong\blinfty{X\cross X,\mu\cross\mu}$, so we may regard
$\ablinfty\whtens\ablinfty$ as a space of measurable
functions on $X\cross X$.  As we shall see in Theorem \ref{theo:schurmult}
{\it infra}, $u=u'$ in $\ablinfty\whtens\ablinfty$ if and only if $u=u'$
locally marginally almost everywhere.

We let 
\begin{equation}
\fB_{\ablinfty}^\sig(\bdop{\abltwo})=\cbrac{T\in\wbdop{\bdop{\abltwo}}:
\begin{matrix} T(M_\vphi aM_\psi)=M_\vphi T(a)M_\psi  \\
\text{for each } a\iin\bdop{\abltwo} \\
\aand\vphi,\psi\iin\ablinfty\end{matrix}}.
\end{equation}

We can now state the main result of this section.  This establishes a framework
which will be useful in the sequel.

\begin{schurmult} \label{theo:schurmult}
If $(X,\mu)$ is a measure space satisfying the
Radon-Nikodym Theorem, then
the following spaces coincide isometrically:

{\bf (i)}   $\schur{X,\mu}$

{\bf (ii)}  $\blinfty{X,\mu}\whtens\blinfty{X,\mu}$

{\bf (iii)}   $\mults\tee{X,\mu}$.

Moreover, the space in (ii) is naturally completely isometrically
isomorphic to each of the following spaces:

{\bf (iv)} $\fB_{\blinfty{X,\mu}}^\sig(\bdop{\bltwo{X,\mu}})$

{\bf (v)}  $\hilfac{\blone{X,\mu}}{\blinfty{X,\mu}}$.

\noindent Finally, if $(X,\mu)$ is $\sig$-finite, then any element in 
$\blinfty{X,\mu}\whtens\blinfty{X,\mu}$ can be represented as a countable sum
of elementary tensors, and the definition of $\schur{X,\mu}$ can be made using
separable Hilbert spaces only.
\end{schurmult}

We note that if $\mu$ is counting measure on $X$, then 
$\schur{X,\mu}=\schur{X}$,
so this result generalises Theorem \ref{theo:classicschurmult}
and Corollary \ref{cor:classicschurmult1}.
Henceforth we will identify $\schur{X,\mu}$ with 
$\blinfty{X,\mu}\whtens\blinfty{X,\mu}$,
and call elements of this space {\it measurable Schur
multipliers}.

In the case that $(X,\mu)$ is $\sig$-finite, U.\ Haagerup appears to obtain
the equivalence of (i), (ii) and (iii), also with countable index sets
in the definition on $\blinfty{X,\mu}\whtens\blinfty{X,\mu}$, and separable
Hilbert spaces in the definition of $\schur{X,\mu}$.  See \cite{haagerup80}
or \cite{walter89}.

\smallskip
\proof  As above we let $\mathrm{L}^p$ denote
$\mathrm{L}^p(X,\mu)$ for $p=1,2\oor\infty$.

(i)=(ii)  This follows the proof of
Corollary \ref{cor:classicschurmult1} almost exactly.
However, if $u=\sum_{i\in I}\vphi_i\tens\psi_i$ in $\ablinfty\whtens\ablinfty$,
we want to ensure that the analogues of (\ref{eq:vectorfromfunc}) are
weakly measurable.  Indeed, let $\xi,\eta:X\to\ltwo{I}$ be given by
\begin{equation}\label{eq:vffromwh}
\xi(x)=\brac{\vphi_i(x)}_{i\in I}\quad\aand\quad
\eta(x)=\brac{\wbar{\psi_i(x)}}_{i\in I}
\end{equation}
for locally almost every $x\iin X$.  Then if $\zeta=\brac{\zeta_i}_{i\in I}$
in $\ltwo{I}$, we have that $\zeta_i=0$ for all but countably many $i$
and hence $\inprod{\xi(\cdot)}{\zeta}=\sum_{i\in I}\wbar{\zeta_i}\vphi_i$
is measurable.  Similarly $\eta$ is weakly measurable.

(ii)$\cong$(iv) The map (\ref{eq:tmap}) restricts to a surjective complete
isometry from $\ablinfty\whtens\ablinfty$ to
$\fC\fB_{\ablinfty}^\sig(\bdop{\abltwo})$ by Corollary 
\ref{cor:ehtensprodchar1} and \cite[Theo.\ 3.1]{smith91}.
Then $\fC\fB_{\ablinfty}^\sig(\bdop{\abltwo})
=\fB_{\ablinfty}^\sig(\bdop{\abltwo})$
by \cite[Theo.\ 3.1]{smith91}, since $\ablinfty$ is locally cyclic
on $\abltwo$.

(ii)$\cong$(v) From (\ref{eq:lonehdual}), 
$\ablinfty\whtens\ablinfty\cong(\ablone\htens\ablone)^*$, completely 
isometrically.
Also, $(\ablone\htens\ablone)^*\cong\hilfac{\ablone}{\ablinfty}$
via (\ref{eq:blonehd}).  We note that the complete isomorphism
$\ablinfty\whtens\ablinfty\cong\hilfac{\ablone}{\ablinfty}$ is
given by $u\mapsto T(u)$, where for $f\iin\ablone$, $T(u)f=L_f u$, where
$L_f$ is the left slice map given in (\ref{eq:slicemap}).  Indeed,
if $u=\sum_{i\in I}\vphi_i\tens\psi_i$,
the factorization $T=T_2\comp T_1$, as in (\ref{eq:hilfac}),
with $\hil=\ltwo{I}$, is given by
\[
T_1f=\brac{\dpair{f}{\vphi_i}}_{i\in I}\quad\aand\quad\
T_2\zeta=\sum_{i\in I}\zeta_i\psi_i
\]
where $\zeta=(\zeta_i)_{i\in I}$.

(ii)=(iii) (via (ii)$\cong$(iv)) Let $u\in\ablinfty\whtens\ablinfty$,
with $u=\sum_{i\in I}\vphi_i\tens\psi_i$ and
$\whnorm{u}=\unorm{\sum_{i\in I}|\vphi_i|^2}^{1/2}
\unorm{\sum_{i\in I}|\psi_i |^2}^{1/2}$.  
It suffices to show that
$m_uf\tens g\in\tee{X,\mu}$ with
$\norm{m_uf\tens g}_\gam\leq\whnorm{u}\norm{f}_2\norm{g}_2$
for any $f,g\iin\abltwo$, to establish that $u\in\mults\tee{X,\mu}$
with $\norm{m_u}\leq\whnorm{u}$.  Using the definition of the projective
tensor product and then H\"{o}lder's Inequality we have that
\begin{align*}
\norm{m_uf\tens g}_\gam
&=\norm{\sum_{i\in I}M_{\psi_i}f\tens M_{\vphi_i}g}_\gam
\leq\sum_{i\in I}\norm{M_{\psi_i}f}_2\norm{M_{\vphi_i}g}_2 \\
&\leq\brac{\sum_{i\in I}\norm{M_{\psi_i}f}_2^2}^{1/2}
\brac{\sum_{i\in I}\norm{M_{\vphi_i}g}_2^2}^{1/2}.
\end{align*}
Then, using Tonelli's Theorem, we have that
\begin{align}
\label{eq:multsum}
\sum_{i\in I}\norm{M_{\psi_i}f}_2^2
&=\sum_{i\in I}\int_X |\psi_i f|^2d\mu
=\int_X \sum_{i\in I}|\psi_i f|^2d\mu \notag \\
&=\int_X \sum_{i\in I}|\psi_i|^2|f|^2d\mu
\leq\unorm{\sum_{i\in I}|\psi_i |^2}\norm{f}_2^2
\end{align}
and, similarly, $\sum_{i\in I}\norm{M_{\vphi_i}g}_2^2\leq
\unorm{\sum_{i\in I}|\vphi_i|^2}\norm{g}_2^2$.  Hence
\[
\norm{m_uf\tens g}_\gam
\leq\unorm{\sum_{i\in I}|\vphi_i|^2}^{1/2}
\unorm{\sum_{i\in I}|\psi_i |^2}^{1/2}\norm{f}_2\norm{g}_2
=\whnorm{u}\norm{f}_2\norm{g}_2
\]
as required.

Now suppose that $u\in\mults\tee{X,\mu}$.  For $\vphi$ and $\psi\iin\ablinfty$,
the adjoint of $M_\psi\tens M_\vphi$, qua operator on $\tee{X,\mu}$, in
the dual pairing described in (\ref{eq:tdualisb}), is $T_{\vphi\tens\psi}$.
Indeed we have for $a\iin\bdop{\abltwo}$ and $f,g\iin\abltwo$,
\begin{align*}
\dpair{a}{M_\psi\tens M_\vphi(f\tens g)}&=\dpair{a}{M_\psi f\tens M_\vphi g}
=\inprod{aM_\psi f}{M_\vphi^*\wbar{g}} \\
&=\inprod{M_\vphi aM_\psi f}{g}
=\dpair{T_{\vphi\tens\psi}a}{f\tens g}.
\end{align*}
Since $M_\psi\tens M_\vphi m_u=m_u M_\psi\tens M_\vphi$ for all 
$\vphi$ and $\psi\iin\ablinfty$, 
$T_{\vphi\tens\psi}m_u^*=m_u^*T_{\vphi\tens\psi}$.
Hence $m_u^*\in\fB_{\ablinfty}^\sig(\bdop{\abltwo})$ and there is
$u'\iin\ablinfty\whtens\ablinfty$ such that $m_u^*=T_{u'}$.  But then
$m_{u'}=m_u$ so $u'=u$ locally marginally almost everywhere, and hence
$u\in\ablinfty\whtens\ablinfty$.  

We finally deal with the case where $(X,\mu)$ is $\sigma$-finite.  We can 
decompose $X=\dot{\bigcup}_{n\in\En}E_n$ where each $E_n$ is measurable with 
$\mu(E_n)<+\infty$.  Then, given $u=\sum_{i\in I}\vphi_i\tens\psi_i$
in $\ablinfty\whtens\ablinfty$, using (\ref{eq:multsum}) we have
\[
\sum_{i\in I}\int_X\abs{\psi_i\chi_{E_n}}^2d\mu
\leq\unorm{\sum_{i\in I}|\psi_i|^2}\mu(E_n)
\]
so $\psi_i\chi_{E_n}=0$ (almost everywhere) for all but countably many indices
$i$.  It then follows that $\psi_i\not=0$ for at most countably many indices
$i$.  The same holds for the functions $\vphi_i$.  
Hence if $u\in\ablinfty\whtens\ablinfty$,
we can write $u=\sum_{i=1}^\infty\vphi_i\tens\psi_i$.  Following
(\ref{eq:vffromwh}), it is clear that only separable Hilbert spaces
are needed in the definition of $\schur{X,\mu}$.  

(The author would
like to express his gratitude to L.\ Turowska for suggesting the proof
of (ii)=(iv), given here.  This is a generalised version
of the proof of \cite[Prop.\ 4.1]{spronkt}.)  \endpf

\smallskip
We observe that, in general, we can bound the size of $I$ by
$|I|\leq\aleph_0\mult c(X,\mu)$, where $\aleph_0$ is the first infinite
cardinal and $c(X,\mu)$ is the smallest cardinality of a cover of
$X$ by $\mu$-finite sets.

We close this section with an asymmetric version of Theorem 
\ref{theo:schurmult}.  If $(X,\mu)$ and $(Y,\nu)$ are measure spaces,
each satisfying the Radon-Nikodym Theorem, then we can consider the
following spaces:

{\bf (i')} $\schur{X, Y; \mu, \nu}$

{\bf (ii')} $\blinfty{X,\mu}\whtens\blinfty{Y,\nu}$

{\bf (iii')} ${_{\blinfty{X,\mu}}}{\mathcal{B}_{\blinfty{Y,\nu}}}
(\bdop{\bltwo{Y,\nu},\bltwo{X,\mu}})$ \quad(left $\blinfty{X,\mu}$ and
right $\blinfty{Y,\nu}$-module maps)

{\bf (iv')} $\hilfac{\blone{X,\mu},\blinfty{Y,\nu}}$

{\bf (v')} $\mults\tee{Y, X ; \nu, \mu}$, where
$\tee{Y,X;\nu,\mu}=\bltwo{Y,\nu}\tens^{\gam}\bltwo{X,\mu}$.

\noindent The space in (i') is defined analogously to that in 
(\ref{eq:mschur});
note that $\schur{X,\mu}=\schur{X,X;\mu,\mu}$.  Then the spaces
(i')-(v') are isometrically isomorphic.  Moreover, the spaces
(ii'), (iii') and (iv') have natural operator space structures in which
they are completely isometrically isomorphic.

That these spaces are isomorphic
can be proved just as Theorem \ref{theo:schurmult}, {\it mutatis 
mutandis}.  Also note that if we take the disjoint union measure space
$(X\dot{\cup}Y,\mu\dot{\cup}\nu)$, then all of the spaces (i')-(v')
can be taken to be ``corners'' of the spaces (i)-(v) of Theorem 
\ref{theo:schurmult}, but with the latter mentioned spaces taken over the 
measure space $(X\dot{\cup}Y,\mu\dot{\cup}\nu)$.

\subsection{Continuous Schur Multipliers}
\label{ssec:contschur}         

For this subsection we will let $X$ and $Y$ denote locally compact Hausdorff
spaces.  Since it will add no extra difficulties in notation, we will
state asymmetrical versions of all our results.

Let
\begin{equation}
\label{eq:schurcocb}
\schurco{X,Y}=\conto{X}\htens\conto{Y}\;\aand\;
\schurcb{X,Y}=\contb{X}\htens\contb{Y}.
\end{equation}
These are clearly closed subalgebras of $\linfty{X}\whtens\linfty{Y}$
and hence Banach algebras under pointwise operations.
If $X$ and $Y$ are compact we will write $\val{X,Y}=\schurco{X,Y}
=\schurcb{X,Y}$.  This algebra is called the {\it Varopoulos algebra}
in \cite{spronkt}.

\begin{cschurmult}
\label{prop:cschurmult}
The Banach algebra $\schurco{X,Y}$ is regular with
Gel'fand spectrum $X\cross Y$.  The same is true for $\schurcb{X,Y}$,
but with spectrum $\sccomp{X}\cross\sccomp{Y}$, where 
$\sccomp{X}\aand\sccomp{Y}$ denote the Stone-\v{C}ech compactifications
of $X$ and $Y$.
\end{cschurmult}

\proof By Grothendieck's Inequality (see \cite[14.5]{defantf}, for example),
$\schurco{X,Y}=$ $\conto{X}\tens^\gam\conto{Y}$.  Since $\conto{X}$ and
$\conto{Y}$ are regular and satisfy the metric approximation property,
the proposition follows from \cite{tomiyama}.  The result for $\schurcb{X,Y}$
follows immediately from the identies $\contb{X}\cong\cont{\sccomp{X}}$
and $\contb{Y}\cong\cont{\sccomp{Y}}$ and the reasoning above. \endpf

\smallskip
In the sequel we will need spaces which extend $\schurco{X,Y}$ and
$\schurcb{X,Y}$.  Let
\begin{equation}
\label{eq:schurob}
\schuro{X,Y}=\conto{X}\ehtens\conto{Y}\;\aand\;
\schurb{X,Y}=\contb{X}\ehtens\contb{Y}.
\end{equation}
It follows from Corollary \ref{cor:ehtensprodchar1} that these
spaces are closed subalgebras of $\linfty{X}\whtens\linfty{Y}$
and hence Banach algebras under pointwise operations.
Moreover, we obtain from Theorem \ref{theo:ehtensprodchar}, the following.

\begin{exthaagerup2}
\label{prop:exthaagerup2}
A function $u:X\cross Y\to\Cee$ defines an element of $\schuro{X,Y}$
(respectively $\schurb{X,Y}$), if and only if there exist families of functions
$\{\vphi_i\}_{i\in I}$ in $\conto{X}$ (respectively $\contb{X}$) and 
$\{\psi_i\}_{i\in I}$ in $\conto{Y}$ (respectively $\contb{Y}$)
such that $u(x,y)=\sum_{i\in I}\vphi_i(x)\psi_i(y)$ for each
$(x,y)\iin X\cross Y$ (i.e. $u=\sum_{i\in I}\vphi_i\tens\psi_i$),
and 
\[
\unorm{\sum_{i\in I}|\vphi_i|^2}\unorm{\sum_{i\in I}|\psi_i|^2}<+\infty.
\]
Moreover the norm of $u$ in $\schuro{X,Y}$ (respectively $\schurb{X,Y}$)
is then
\[
\vnorm{u}=\min\cbrac{\unorm{\sum_{i\in I}|\vphi_i|^2}^{1/2}
\unorm{\sum_{i\in I}|\psi_i|^2}^{1/2}:u=\sum_{i\in I}\vphi_i\tens\psi_i}.
\]
\end{exthaagerup2}

\smallskip
Observe that we can obtain a bound for
the cardinality of the index set $I$:  
$|I|\leq\min\{d(\contb{X}),d(\contb{Y})\}$, where for any Banach space
$\fX$, $d(\fX)$ denotes the {\it density degree} of $\fX$.
If $X$ and $Y$ admit a $\sig$-finite regular Borel measures $\mu$ and $\nu$
such that $\mu(U)>0$ and $\nu(V)>0$ for each open subset $U$ of $X$ and $V$
of $Y$, then $\schurb{X,Y}$ (and hence $\schuro{X,Y}$) is a closed
subspace of $\schur{X,Y;\mu,\nu}$, and hence $I$ can be taken to be countable
by Theorem \ref{theo:schurmult}.

We will find the following description of $\schurb{X,Y}$ useful.
Let for a Hilbert space $\hil$,
\begin{equation*}
\contb{X,\weaks\hil}=\cbrac{\xi:X\to\hil:
\begin{matrix} \xi\text{ is continuous when }\hil\text{ has the weak*} \\
\text{ topology, and }\unorm{\xi}=\sup_{x\in X}|\xi(x)|<+\infty\end{matrix}}.
\end{equation*}
Similarly, define $\contb{Y,\weaks\hil}$.  If $\xi\in\contb{X,\weaks\hil}$
and $\eta\in\contb{Y,\weaks\hil}$, let $u_{\xi,\eta}:X\cross Y\to\Cee$ be
given by $u_{\xi,\eta}(x,y)=\inprod{\xi(x)}{\eta(y)}$.

\begin{exthaagerup3}
\label{prop:exthaagerup3}
If $X$ and $Y$ are locally compact Hausdorff spaces, then
$u\in\schurb{X,Y}$ if and only if there are a Hilbert space $\hil$ and 
functions $\xi\iin\contb{X,\weaks\hil}$ and $\eta\iin\contb{Y,\weaks\hil}$
such that $u=u_{\xi,\eta}$.
Moreover, $\vnorm{u}=\min\cbrac{\unorm{\xi}\unorm{\eta}:
u=u_{\xi,\eta}\text{ as above}}$.
\end{exthaagerup3}

\proof This follows the proof of Corollary \ref{cor:classicschurmult1}
almost exactly.  
However, if $u=\sum_{i\in I}\vphi_i\tens\psi_i$ in $\schurb{X,Y}$,
we want to ensure that the analogues of (\ref{eq:vectorfromfunc}) are
suitably continuous.  Indeed, let $\xi:X\to\ltwo{I}$ and $\eta:Y\to\ltwo{I}$
be given by
\[
\xi(x)=\brac{\vphi_i(x)}_{i\in I}\quad\aand\quad
\eta(y)=\brac{\wbar{\psi_i(y)}}_{i\in I}
\]
for each $x\iin X$ and $y\iin Y$.
To see that $\xi\in\contb{X,\weaks\hil}$, let
$\zeta=(\zeta_i)_{i\in I}$ in $\ltwo{I}$.  
Then we have that
\[
\abs{\inprod{\xi(\cdot)}{\zeta}}=\abs{\sum_{i\in I}\zeta_i\vphi_i}
\leq\norm{\zeta}_2\unorm{\sum_{i\in I}|\vphi_i|^2}^{1/2}
\]
so the middle series converges uniformly, and hence 
$\inprod{\xi(\cdot)}{\zeta}$ is continuous.
Similarly, $\eta\in\contb{Y,\weaks\hil}$. \endpf

\smallskip
Although an element $u\iin\schurb{X,Y}$ always unambiguously describes
a function on $X\cross Y$, it may not describe a continuous function.
This fact is pointed out in \cite{shulmant}.
As an example, let us take $X=Y=[0,1]$, the unit interval in $\Ree$.  For each
$n\iin\En$ let $\vphi_n=\psi_n:[0,1]\to[0,1]$ be a continuous function 
supported on $\sbrac{\frac{1}{2^{n+1}},\frac{1}{2^n}}$ such that 
$\vphi_n\brac{\frac{3}{2^{n+2}}}=\psi_n\brac{\frac{3}{2^{n+2}}}=1$ 
(where $\frac{3}{2^{n+2}}$ is simply the 
midpoint of $\frac{1}{2^{n+1}}$ and $\frac{1}{2^n}$).  Then it is clear that
\[
\unorm{\sum_{n=1}^\infty|\vphi_n|^2}\unorm{\sum_{n=1}^\infty|\psi_n|^2}
=1\cdot 1=1
\]
and hence that $u=\sum_{n=1}^\infty \vphi_n\tens\psi_n\in\schurb{[0,1],[0,1]}$
with norm $1$.  However, while $u(0,0)=0$, 
$u\brac{\frac{3}{2^{n+2}},\frac{3}{2^{n+2}}}=1$
for each $n$, and hence $u$ is not continuous at $0$.  We note, then, that
for any choice of spaces $X\aand Y$,
finding the Gel'fand spectrum for $\schurb{X,Y}$ or $\schuro{X,Y}$, should
prove to be challenging.

It would be interesting to have intrinsic
descriptions of $\schurb{X,Y}\cap\contb{X\cross Y}$,
$\schuro{X,Y}\cap\contb{X\cross Y}$ and $\schuro{X,Y}\cap\conto{X\cross Y}$.

\section{Completely Bounded Multipliers of the Fourier Algebra}
\label{sec:cbmults}

For this section the symbol $\fH$, possibly with subscript,
will denote a Hilbert space, and $\fU(\hil)$ the group of unitaries on $\hil$,
topologized with the weak operator topology.

Our standard references for harmonic analysis are \cite{hewittrI}, 
\cite{eymard} and \cite{arsac}.

Let $G$ be a locally compact group with left Haar integral
$f\mapsto\int_G f(s)ds$.  For $p=1,2$, let $\mathrm{L}^p(G)$ denote
the space of (almost everywhere equivalence classes of) functions
$f$ for which $\norm{f}_p=\brac{\int_G |f(s)|^pds}^{1/p}<+\infty$.
Then $\bloneg$ is an involutive Banach algebra under the product of 
convolution, $f\con g(t)=\int_G f(s)g(s^{-1}t)ds$ for almost every $t$,
and with involution $f^*(t)=\frac{1}{\Del(t)}\wbar{f(t^{-1})}$, where
the overline denotes complex conjugation and $\Del:G\to\Ree^{>0}$ is
the Haar modular function.

If $\pi$ is a continuous unitary representation of $G$, i.e.\ a continuous
homomorphism from $G$ into $\fU(\hil_\pi)$, then it forms a canonical 
non-degenerate
involutive representation $\pi_1:\bloneg\to\bdop{\hil_\pi}$ via the
weak*-converging integrals $\pi_1(f)=\int_Gf(s)\pi(s)ds$.  
We let $\cstaro{\pi}$ be the norm closure of $\pi_1(\bloneg)$
in $\bdop{\hil_\pi}$, and
$\cstarg$ denote the {\it enveloping C*-algebra} of $G$, i.e.\ the enveloping 
C*-algebra of $\bloneg$.

The {\it Fourier} and {\it 
Fourier-Stieltjes algebras}, $\falg$ and $\fsalg$, are defined in 
\cite{eymard}. We recall that $\fsalg$ is the space of matrix coefficients 
of all continuous unitary representations of $G$; i.e.\ the space of functions 
of the form $s\mapsto\inprod{\pi(s)\xi}{\eta}$ where $\pi$ is a continuous
unitary representation of $G$ on $\hil_\pi$, and $\xi,\eta\in\hil_\pi$.
$\fsalg$ is the dual of $\cstarg$ via 
$\dpair{a}{\inprod{\pi(\cdot)\xi}{\eta}}=\inprod{\pi_*(a)\xi}{\eta}$, where
$\pi_*:\cstarg\to\bdop{\fH_\pi}$ is the representation induced by $\pi$, 
i.e.\ by $\pi_1$.  We will refer to $\sig(\fsalg,\cstarg)$ as {\it the}
weak* topology on $\fsalg$, although its uniqueness is not assured.
It is well-known that $\fsalg$ is a commutative Banach algebra (under pointwise
operations) and $\falg$ is the closed ideal in $\fsalg$ generated by
compactly supported matrix coefficients.

If $\pi$ is a continuous unitary representation of $G$, let $\ccfs{\pi}$
be the norm closure of 
$\spn\{\inprod{\pi(\cdot)\xi}{\eta}:\xi,\eta\in\fH_\pi\}$ in $\fsalg$.
Then, by \cite[2.2]{arsac}, $\ccfs{\pi}^*\cong\vn{\pi}$, where $\vn{\pi}$ is 
the von Neumann algebra generated by $\pi(G)$.  Thus if $\vpi$ denotes
the universal representation of $G$, we have $\fsalg=\ccfs{\vpi}$.
The algebra $\wstarg=\vn{\vpi}$ is called the {\it enveloping von Neumann
algebra} of $G$.  Also, $\falg=\ccfs{\lam}$ where $\lam$ is the left
regular representation of $G$.  We let $\vng=\vn{\lam}$, and call it the
{\it group von Neumann algebra}.  We will also let $\wcfs{\pi}$ denote
the weak*-closure of $\ccfs{\pi}$ in $\fsalg$.  By standard functional
anaysis, $\wcfs{\pi}\cong\bigl(\cstaro{\pi}\bigr)^*$.  In particular, we let
$\fsalrg=\wcfs{\lam}$.  Then $\fsalrg$ is an ideal in $\fsalg$, and is
the dual of the {\it reduced group C*-algebra} $\cstarrg=\cstaro{\lam}$.

\subsection{Multipliers and Completely Bounded Multipliers}

A {\it completely contractive Banach algebra} is a Banach algebra $\fA$, 
equipped with an operator space structure such that the multiplication map
$\m:\fA\pptens\fA\to\fA$ is a complete contraction, and hence extends
to a complete contraction $\m:\fA\ptens\fA\to\fA$.

The Fourier-Stieltjes algebra $\fsalg$, as the dual of $\cstarg$, 
is a completely contractive Banach 
algebra.  Indeed, the multiplication map $\m:\fsalg\pptens\fsalg\to\fsalg$
has adjoint $\m^*:\wstarg\to\wstarg\vntens\wstarg$, given for $s\iin G$ by
$\m^*\vpi(s)=\vpi(s)\tens\vpi(s)$.  By the universal property of
$\wstarg$, $\m^*$ is a $*$-homomorphism and hence a complete contraction;
thus so too is $\m$.

Any of the subspaces 
$\ccfs{\pi}$ of $\fsalg$ obtains the same operator space structure
qua subspace of $\fsalg$ as ii does as the predual of $\vn{\pi}$.  
Hence if $\ccfs{\pi}$ is a subalgebra of $\fsalg$, then it is 
a completely contractive Banach algebra.  In particular, $\falg$
is a completely contractive Banach algebra.

Let $\fA$ be a semi-simple commutative Banach algebra with spectrum $X$, or
more generally, a subalgebra of $\contb{X}$ for which there is some
$a\iin\fA$ such that $a(x)\not=0$ for each $x\iin X$, and comes
equipped with a norm under which it is a Banach algebra.  The {\it multipliers}
of $\fA$ are given by
\[
\mults\fA=\{u:X\to\Cee:ua\in\fA\text{ for each }a\iin\fA\}.
\]
An application of the Closed Graph Theorem shows that the map $m_u:\fA\to\fA$,
given by 
\begin{equation}
\label{eq:multu}
m_u a=ua
\end{equation}
is bounded.  Moreover, it can be shown that $\{m_u:u\in\mults\fA\}$
is the strong operator topology closed subalgebra of $\fA$-module maps
in $\bdop{\fA}$.  See \cite{larsenB}, for example.  If $\fA$ has an
operator space structure let
\[
\cbmults\fA=\{u\in\mults\fA:m_u\in\cbop{\fA}\}.
\]
The operator space structure is given as follows:
for $\sbrac{u_{ij}}\iin
\matn{\cbmults\fA}$ let
\[
\cbmnorm{\sbrac{u_{ij}}}=\cbnorm{m_{\sbrac{u_{ij}}}}
\]
where $m_{\sbrac{u_{ij}}}\iin\cbop{\fA,\matn{\fA}}$ is given by
\begin{equation}
m_{\sbrac{u_{ij}}}a=\sbrac{m_{u_{ij}}a}=\sbrac{u_{ij}a}
\end{equation}
for $a\iin\fA$.

If $\fA$ is, in addition to the properties assumed above, 
a completely contractive Banach algebra, then it injects
completely contractively into $\cbmults\fA$.  Moreover, if 
$\fA$ has a contractive bounded approximate identity $\{e_\alp\}$, then
the injection $\fA\hookrightarrow\cbmults\fA$ is a complete isometry.
Indeed, if $\sbrac{a_{ij}}\in\matn{\fA}$ then
\begin{equation}
\label{eq:comisominj}
\norm{\sbrac{a_{ij}}}\geq\cbmnorm{\sbrac{a_{ij}}}\geq
\lim_\alp\norm{\sbrac{a_{ij}e_\alp}}=\norm{\sbrac{a_{ij}}}.
\end{equation}

Now let $G$ be a locally compact group and $\falg$ the Fourier algebra
of $G$.  The following Proposition can be proved exactly as 
\cite[Prop. 1.2]{decanniereh}, but using the fact that the adjoint
map $T\mapsto T^*$ from $\cbop{\fV}$ to $\wcbop{\fV^*}$ is a complete
isometry for any operator space $\fV$, and that the restriction map
$T\mapsto T|_{\fV_0}$ form $\cbop{\fV}$ to $\cbop{\fV_0,\fV}$, for any
closed subspace $\fV_0$, is a complete contraction 

\begin{fmultschar}
\label{prop:fmultschar}
If $u_{ij}:G\to\Cee$ for $i,j=1,\dots,n$, then the following are equivalent:

{\bf (i)}  $\sbrac{u_{ij}}\in\ball{\matn{\cbmg}}$.

{\bf (ii)} The operator defined
by $\lam(s)\mapsto\sbrac{u_{ij}(s)\lam(s)}$, for $s\iin G$, extends to
a weakly continuous complete contraction, $M_{\sbrac{u_{ij}}}:\vng
\to\matn{\vng}$.

{\bf (iii)} The operator defined by  
$\lam_1(f)\mapsto\sbrac{\lam_1(u_{ij}f)}$,
for $f\iin\bloneg$, extends to a complete contraction 
$\til{M}_{\sbrac{u_{ij}}}:\cstarrg\to\matn{\cstarrg}$.

{\bf (iv)} $\sbrac{u_{ij}}\in\ball{\matn{\cbmults\fsalrg}}$.

\noindent In (ii), $M_{\sbrac{u_{ij}}}T=\sbrac{M_{u_{ij}}T}$, for $T\iin\vng$, 
where for $i,j=1,\dots,n$, $M_{u_{ij}}\in\wcbop{\vng}$
is the adjoint of $m_{u_{ij}}\iin\cbop{\falg}$.
\end{fmultschar}

Following (iii) of the above proposition, we may consider $\cbmg$ to be
the ``multipliers of $\cstarrg$''.  We will see below that 
$\fsalg$ represents the multipliers of the enveloping C*-algebra,
$\cstarg$.  In \cite{walter86} and \cite{walter89}, 
it is shown that $\fsalg$ imbeds as a maximal Abelian
subalgebra of $\cbop{\cstarg}$.  We prove this below, and also show that
the operator space structure, induced on $\fsalg$ as a subspace of
$\cbop{\cstarg}$, is the same as the operator space structure induced on 
$\fsalg$ as the dual of $\cstarg$.

\begin{fsalgchar}
\label{prop:fsalgchar}
If $u_{ij}:G\to\Cee$ for $i,j=1,\dots,n$, then the following are equivalent:

{\bf (i)}  $\sbrac{u_{ij}}\in\ball{\matn{\fsalg}}$.
 
{\bf (ii)} $\sbrac{u_{ij}}\in\ball{\matn{\cbmults\fsalg}}$.

{\bf (iii)} The operator defined
by $\vpi(s)\mapsto\sbrac{u_{ij}(s)\vpi(s)}$, for $s\iin G$, extends to
a weakly continuous complete contraction
$M_{\sbrac{u_{ij}}}:\wstarg\to\matn{\wstarg}$.

{\bf (iv)} The operator defined by  $\vpi_1(f)\mapsto\sbrac{\vpi_1(u_{ij}f)}$,
for $f\iin\bloneg$, extends to a complete contraction
$\til{M}_{\sbrac{u_{ij}}}:\cstarg\to\matn{\cstarg}$.

\noindent Thus $\fsalg$ is completely isometrically isomorphic to the maximal 
Abelian subalgebra $\{\til{M}_u:u\in\fsalg\}$ of $\cbop{\cstarg}$.
\end{fsalgchar}

\proof The fact that (i) is equivalent to (ii) follows (\ref{eq:comisominj}).
That (ii), (iii) and (iv) are equivalent can be shown exactly as
the equivalence of (ii)-(iv) in the preceding proposition, i.e.\ by
following \cite[Prop. 1.2]{decanniereh}.

We have that $\{\til{M}_u:u\in\fsalg\}$ is a maximal Abelian subalgebra 
of $\cbop{\cstarg}$, since if $T\iin\cbop{\cstarg}$ commutes with every
$\til{M}_u$, where $u\in\fsalg$, then $T^*$ is a $\fsalg$-module map on 
$\fsalg$, and hence $T=m_u$ for some $u\iin\fsalg$. \endpf

\smallskip
We note that, in fact, $\{\til{M}_u:u\in\fsalg\}$ is a maximal Abelian
subalgebra of $\bdop{\cstarg}$.

We can refine the well known result that $\fsalg\subset\cbmg$ in general, and
that $\fsalg=\cbmg$ when $G$ is amenable.
See \cite[Cor.\ 1.8]{decanniereh} and \cite{losert84}, for example.

\begin{fsalgchar1}\label{cor:fsalgchar1}
$\fsalg\subset\cbmg$, and the identity injection $\fsalg\hookrightarrow\cbmg$ 
is a complete contraction.  If $G$ is amenable, then $\iota$ is a completely 
isometric isomorphism.
\end{fsalgchar1}

\proof  If $u\in\fsalg$ let $m_u$ be the multiplier operator in 
$\cbop{\fsalg}$.  Then, $m_u\falg\subset\falg$, since $\falg$ is an ideal in
$\fsalg$.  Moreover, $m_u\mapsto m_u|_{\falg}$ is a complete contraction.

If $G$ is amenable, then $\fsalg=\fsalrg$, completely isometrically.
the result then follows from Proposition \ref{prop:fmultschar} (iv) 
and Proposition \ref{prop:fsalgchar} (ii). \endpf

\section{Invariant Schur Multipliers}
\label{sec:invschur}

For this section, let $G$ be a locally compact group with
left Haar measure $m$.  Let $\schur{G}=\schur{G,m}$.  Hence, if $G$ is 
discrete, $\schur{G}$ is just the space of usual Schur multipliers on
$\bdop{\ell^2(G)}$.

Let
\begin{equation}
\schurinv{G}=\cbrac{u\in\schur{G}:\begin{matrix}
\text{for all }r\iin G, u(sr,t)=u(s,tr^{-1}) \\
\text{for l.m.a.e. }(s,t)\iin G\cross G \end{matrix}}
\end{equation}
In the identification $\schur{G}=\blinftyg\whtens\blinftyg\cong
(\bloneg\htens\bloneg)^*$ in Theorem \ref{theo:schurmult} we have that
\begin{equation*}
\schurinv{G}=\bigcap_{r\in G}
\ker\brac{\id_{\blinftyg}\tens\rho(r^{-1})^*-\rho(r)^*\tens\id_{\blinftyg}}.
\end{equation*}
Here $\rho:G\to\bdop{\blone{G}}$ is the right regular representation  
on $\blone{G}$, given by $\rho(r)f(s)=\Del(r)f(sr)$, for
$r\iin G$, $f\iin\blone{G}$ and almost all $s\iin G$.
Hence $\schurinv{G}$ is a weak* closed subspace of $\schur{G}$.

If $G$ is a discrete group, then $\schur{G}$ is a space of functions, rather
that a space of equivalence classes of functions.  Thus
for $u\iin\schurinv{G}$ it is evident
that we can define a function $u_G:G\to\Cee$ such that
\[
u_G(st^{-1})=u(s,t)
\]
for all $(s,t)\iin G\cross G$.  For a non-discrete group this is more
subtle.

\begin{invschurmult1}\label{lem:invschurmult1}
Suppose $G$ is a non-discrete locally compact group, and
let $u\in\schurinv{G}$.  Consider $u$ as a function on $G\cross G$, i.e.\
as a representative of its equivalence class.
Then there is a locally null subset $E_u$ of $G$ such that
for $t\iin G\setdif E_u$, $u(s,t^{-1}s)=
u(e,t^{-1})$ for locally almost every $s\iin G$. 
\end{invschurmult1}

\proof Let for $t\iin G$, $F_t=\{s\in G:u(s,t^{-1}s)\not
=u(e,t^{-1})\}$.  Then let $E_u=\{t\in G:F_t\text{ is not locally null}\}$.
We will show that $E_u$ is locally null.
If $t\in E_u$ and $s\in F_t$, then $u(s,t^{-1}s)\not=u(e,t^{-1})
=u(ss^{-1},t^{-1})$, so there is an $r\iin G$, namely $r=s^{-1}$, such
that $u(sr,t^{-1})\not=u(s,t^{-1}r^{-1})$.  Hence
\[
\bigcup_{t\in E_u} F_t\cross\{t^{-1}\}\subset
\{(s,t)\in G\cross G:u(sr,t^{-1})\not=u(s,t^{-1}r^{-1})\text{ for some }
r\iin G\}
\]
and thus $\bigcup_{t\in E_u} F_t\cross\{t^{-1}\}$ is locally marginally
null.  Let $N_1$ and $N_2$ be locally null subsets of $G$ such that
\[
\bigcup_{t\in E_u} F_t\cross\{t^{-1}\}\subset(N_1\cross G)\cup(G\cross N_2).
\]
Note that $F_t\setdif N_1\not=\varnothing$ for each $t\iin E_u$, since,
for such $t$, $F_t$ is not locally null.   Then, since
$\brac{\bigcup_{t\in E_u} (F_t\setdif N_1)\cross\{t^{-1}\}}\cap(N_1\cross G)
=\varnothing$, we have that
\[
\bigcup_{t\in E_u} (F_t\setdif N_1)\cross\{t^{-1}\}\subset G\cross N_2.
\]
Hence $E_u^{-1}\subset N_2$ so $E_u\subset N_2^{-1}$, and thus $E_u$ is locally
null.  \endpf

\smallskip
If $u\in\schurinv{G}$, considered as a function as in the above lemma, 
but bounded, define for $t\iin G$,
\begin{equation}\label{eq:usubg}
u_G(t)=\begin{cases}
u(e,t^{-1}) &\iif t\in G\setdif E_u \\
\underset{s\in V,V\in\fV_t}{\lim}u(e,s^{-1}) &\iif t\in E_u
\end{cases}
\end{equation}
where for $t\iin E_u$, $\fV_t$ is any fixed ultra-filter over $G\setdif E_u$ 
containing the sets $U\setdif E_u$, where $U$ is a neighbourhood of $t$.  Note
that $U\setdif E_u$ is always non-empty, since an $m$-locally null set in $G$ 
never contains an open set.  It is shown in the next theorem that $u_G$ is
dependent only on the choice of $u\iin\schurinv{G}$ qua equivalence class.
Moreover, it shown that $u_G$ is automatically continuous.

Let $\tee{G}=\bltwog\tens^\gam\bltwog$, as in Subsection 
\ref{ssec:measschurmult}.
Also recall the convention that for $\ome=f\tens g\iin\tee{G}$,
$\ome(s,t)=f(t)g(s)$ for m.a.e.\ $(s,t)\iin G\cross G$.

\begin{invschurmult2}\label{theo:invschurmult2}
For each $u\iin\schurinv{G}$, $u_G\in\cbmg$.  
The map $u\mapsto u_G$ induces a complete contraction from 
$\schurinv{G}$ to $\cbmg$.
\end{invschurmult2}

\proof The map $P:\tee{G}\to\falg$ given by
\[
Pf\tens g(t)=\inprod{\lam(t)f}{\wbar{g}}=\int_G f(t^{-1}s)g(s)ds
\]
for $t\iin G$ and $f,g\iin\bltwog$, is a surjective contraction whose
adjoint, $P^*$, is the usual imbedding $\vng\hookrightarrow\bdop{\bltwog}$.
Thus if $u\in\schurinv{G}$, letting $m_u$ be as in (\ref{eq:multstee}),
we have that $P(m_u\ome)\in\falg$ for $\ome\iin\tee{G}$.  If $G$ is not
discrete, let $E_u$ be as the above lemma; if $G$ is discrete let 
$E_u=\varnothing$.  Then for $t\iin G\setdif E_u$  we have
\[
P(m_u\ome)(t)=\int_G u(s,t^{-1}s)\ome(s,t^{-1}s)ds
=u_G(t)P\ome(t).
\]
Let $v\in\falg$ with $v=P\ome$, for some $\ome\iin\tee{G}$.
Since $v$ is continuous, using the notation of (\ref{eq:usubg}), we have that 
\[
u_G(t)v(t)=\lim_{s\in V,V\in\fV_t}u_G(s)v(s)
=\lim_{s\in V,V\in\fV_t}P(m_u\ome)(s)=P(m_u\ome)(t)
\] 
for $t\iin E_u$.  Hence $u_G(t)v(t)=P(m_u\ome)(t)$ for all $t\iin G$, 
so $u_Gv\in\falg$.  Then, by definition, $u_G\in\mults\falg$ with 
$P\comp m_u=m_{u_G}\comp P$, where $m_{u_G}$ in $\bdop{\falg}$ is as defined in
(\ref{eq:multu}).  Taking adjoints we have that $T_u\comp P^*=(P\comp m_u)^*=
(m_{u_G}\comp P)^*=P^*\comp M_{u_G}$, where $M_{u_G}\iin\wbdop{\vng}$ is the 
adjoint of $m_{u_G}$, and hence for $x\iin\vng$ we have that,
\[
T_u x=m_u^*x=M_{u_G}x.
\]
Thus $M_{u_G}=T_u|_{\vng}$, so
$\cbmnorm{u_G}=\cbnorm{M_{u_G}}\leq\cbnorm{T_u}=\scnorm{u}$.
It also follows that for $\sbrac{u_{ij}}\iin\matn{\schurinv{G}}$,
$\cbmnorm{\sbrac{u_{ij}}}$ $\leq\scnorm{\sbrac{(u_{ij})_G}}$. \endpf

\smallskip
The following theorem refines the main result in \cite{bozejkof}. 
We consider this to be the centrepiece of this paper. 
Our proof uses modern operator space theoretic techniques, and
gives us a complete isometry.
Also, it exposes the relationship between $\cbmg$
and $\blinftyg$-module maps on $\bdop{\bltwog}$ via Theorem
\ref{theo:schurmult}.  Many of our methods are motivated
by the elegant exposition found in \cite{haagerup}.

\begin{haagerup}
\label{theo:haagerup}
If $u_{ij}:G\to\Cee$ are continuous functions for $i,j=1,\dots,n$, then
the following are equivalent:

{\bf (i)} $\sbrac{u_{ij}}\in\ball{\matn{\cbmg}}$.

{\bf (ii)} For any finite subset $F$ of $G$, if 
$u_{ij}^F(s,t)=u_{ij}(st^{-1})$
for $s,t\iin F$, then $\sbrac{u_{ij}^F}\in\ball{\matn{\schur{F}}}$.

{\bf (iii)} If $u_{ij}^G$ is given by $u_{ij}^G(s,t)=u_{ij}(st^{-1})$
for $s,t\iin G$, then $\sbrac{u_{ij}^G}\in\ball{\matn{\schurinv{G_d}}}$.

{\bf (iv)} $\sbrac{u_{ij}^G}\in\ball{\matn{\schurinv{G}}}$.

\noindent In particular, the map $u\mapsto u^G$ is a complete isometry from
$\cbmg$ onto $\schurinv{G}$.
\end{haagerup}

\proof (i)$\implies$(ii) Let $F=\{s_1,\dots,s_m\}$, and for $k,l=1,\dots,m$ 
let $E_{kl}$ in $\mat{m}{\vng}$ be identified by
\[
E_{kl}\cong\lam(s_ks_l^{-1})\tens e_{kl}\quad\iin\quad\vng\stens\smat{m}
\cong\mat{m}{\vng}
\]
where $\{e_{kl}\}_{k,l=1,\dots,m}$ is the standard matrix unit for
$\smat{m}$.  Then $\{E_{kl}\}_{k,l=1,\dots,m}$ is itself a matrix unit so
$\fN_F=\spn\{E_{kl}\}_{k,l=1,\dots,m}\cong\smat{m}$ via
$E_{kl}\mapsto e_{kl}$.  Let us identify $\smat{m}\cong\bdop{\ltwo{F}}$.
Note that for $v\iin\schur{F}$, $S_ve_{kl}=v(s_k,s_l)e_{kl}$, for
each $k,l$, where $S_v$ is defined as in (\ref{eq:schurmults}).
If $\sbrac{v_{ij}}\in\matn{\schur{F}}$, let $S_{\sbrac{v_{ij}}}x=
\sbrac{S_{v_{ij}}x}$ for $x\iin\bdop{\ltwo{F}}$.  We then have
\begin{align*}
\lift{m}{M_{\sbrac{u_{ij}}}}E_{kl}
&\cong\sbrac{M_{u_{ij}}\lam(s_ks_l^{-1})\tens e_{kl}}
=\sbrac{u_{ij}(s_ks_l^{-1})\lam(s_ks_l^{-1})\tens e_{kl}} \\
&\cong\sbrac{u_{ij}(s_ks_l^{-1})E_{kl}}
\cong\sbrac{u_{ij}^F(s_k,s_l)e_{kl}}
\cong S_{\sbrac{u_{ij}^F}}e_{kl}
\end{align*}
so $\sbrac{u_{ij}^F}\in\matn{\schur{F}}$.  Moreover,
\[
\scnorm{\sbrac{u_{ij}^F}}=\cbnorm{S_{\sbrac{u_{ij}^F}}}
=\cbnorm{M_{\sbrac{u_{ij}}}|_{\fN_F}}\leq\cbmnorm{\sbrac{u_{ij}}}\leq 1
\]
as required.

(ii)$\implies$(iii) Let $\fF$ denote the collection of all finite subsets of 
$G$, directed by inclusion.  
We have for each $F\iin\fF$, a natural completely isometric imbedding
\[
\schur{F}=\linfty{F}\htens\linfty{F}\subset
\linfty{G_d}\whtens\linfty{G_d}=\schur{G_d}
\]
and hence an isometric imbedding $\matn{\schur{F}}\subset\matn{\schur{G_d}}$.
Denote the image of $\sbrac{u_{ij}^F}$ under this imbedding
again by $\sbrac{u_{ij}^F}$.
Since 
$\matn{\schur{G_d}}$ 
is a dual space (with predual $\tmatn{\lone{G_d}\htens\lone{G_d}}$), the
net $\brac{\sbrac{u_{ij}^F}}_{F\in\fF}$ must have a weak* cluster point
in $\ball{\matn{\schur{G_d}}}$.  Using the fact that
$\bigcup_{F\in\fF}\lone{F}\htens\lone{F}$ is dense in 
$\lone{G_d}\htens\lone{G_d}$, it is easy to see that
\[
\sbrac{u_{ij}^G}=\text{w*}-\lim_{F\in\fF}\sbrac{u_{ij}^F}
\]
so $\sbrac{u_{ij}^G}\in\ball{\matn{\schur{G_d}}}$.  That
$\sbrac{u_{ij}^G}\in\matn{\schurinv{G_d}}$ is clear.

(iii)$\implies$(iv) By (\ref{eq:matrixwhnorm}) we may write
$\sbrac{u_{ij}^G}$ as a matricial product $A\odot B$, where
$A\in\mat{n,I}{\linfty{G_d}}$ and $B\in\mat{I,n}{\linfty{G_d}}$, for some
index set $I$, with $\norm{A}\norm{B}\leq 1$.  We may explicitly write
this down:  well-order $I=\{\iota_1,\iota_2,\dots\}$ so that we have
\begin{equation}\label{eq:matdecomp}
\sbrac{u_{ij}^G}=
\begin{bmatrix}
\vphi_{1\iota_1} & \vphi_{1\iota_2} &\vphi_{1\iota_3} &\cdots \\
\vdots & \vdots & \vdots & \\
\vphi_{n\iota_1} & \vphi_{n\iota_2} &\vphi_{n\iota_3} & \cdots
\end{bmatrix}
\odot
\begin{bmatrix}
\psi_{\iota_1 1} & \cdots & \psi_{\iota_1 n} \\
\psi_{\iota_2 1} & \cdots & \psi_{\iota_2 n} \\
\psi_{\iota_3 1} & \cdots & \psi_{\iota_3 n} \\
\vdots & &\vdots
\end{bmatrix}
\end{equation}
and we see that
\[
u_{ij}^G=\sum_{\iota\in I}\vphi_{i\iota}\tens\psi_{\iota j}.
\]
For $i,j=1,\dots,n$ let $\xi_i,\eta_j:G\to\ltwo{I}$ be given by
\[
\xi_i(s)=(\vphi_{i\iota}(s))_{\iota\in I}\quad\aand\quad
\eta_j(s)=\brac{\wbar{\psi_{\iota j}(s)}}_{\iota\in I}
\]
so that $u_{ij}^G(s,t)=\inprod{\xi_i(s)}{\eta_j(t)}$ for $s,t\iin G$.
If we consider each $\xi_i$ and $\eta_j$ as row vectors with $\linfty{G_d}$
entries, i.e.\ as elements of $\mat{1,I}{\linfty{G_d}}$,
then we may rewrite (\ref{eq:matdecomp}) as
\begin{equation}\label{eq:matdecomp1}
\sbrac{u_{ij}^G}=
\begin{bmatrix} \xi_1 \\ \vdots \\ \xi_n \end{bmatrix}
\odot
\begin{bmatrix} \eta_1^* & \cdots & \eta_n^* \end{bmatrix}.
\end{equation}

Let $\hil'=\wbar{\spn}\{\eta_j(t):j=1,\dots,n,\; t\in G\}$ and
$p'\iin\bdop{\ltwo{I}}$ denote the orthogonal projection onto $\hil'$.
Then let $\hil=\wbar{\spn}\{p'\xi_i(s):i=1,\dots,n,\; s\in G\}$ and
$p\iin\bdop{\ltwo{I}}$ denote the orthogonal projection onto $\hil$.
Clearly $pp'=p$, and $\inprod{p\xi_i(s)}{\eta_i(t)}=
\inprod{\xi_i(s)}{\eta_i(t)}$ for all indices $i$ and $j$, and all 
$s,t\iin G$.  Let $\til{\xi}_i,\til{\eta}_j:G\to\ltwo{I}$ be given by
\[
\til{\xi}_i(s)=p\xi_i(s)\quad\aand\quad\til{\eta}_j(s)=p\eta_j(s).
\]
For any $i=1,\dots,n$ the function $\til{\xi}_i$ is weakly continuous.  
Indeed for
any $s\iin G$ and $\eta\iin\ltwo{I}$, $\inprod{\til{\xi}_i(s)}{\eta}=
\inprod{\xi_i(s)}{p\eta}$ where $p\eta\in\hil'$.  Then, for any
$t\iin G$ and $j=1,\dots,n$, $\inprod{\xi_i(s)}{\eta_j(t)}=u_{ij}^G(s,t)
=u_{ij}(st^{-1})$, where $s\mapsto u_{ij}(st^{-1})$ is continuous 
by assumption.
Similarly, each $\til{\eta}_j$ is continuous.  Thus for any $i, j=1,\dots,n$
and $\iota\iin I$, the functions
\begin{equation}
\label{eq:tilvphitilpsi}
\til{\vphi}_{i\iota}=\inprod{\til{\xi}_i(\cdot)}{e_\iota}\quad
\aand\quad\til{\psi}_{\iota j}=\inprod{e_\iota}{\til{\eta}_i(\cdot)}
\end{equation}
are continuous.  Then, following the notation of (\ref{eq:matdecomp}), if

\[
\til{A}=\begin{bmatrix}
\til{\vphi}_{1\iota_1} & \til{\vphi}_{1\iota_2} 
&\til{\vphi}_{1\iota_3} &\cdots \\
\vdots & \vdots & \vdots & \\
\til{\vphi}_{n\iota_1} & \til{\vphi}_{n\iota_2} 
&\til{\vphi}_{n\iota_3} & \cdots
\end{bmatrix}\aand
\til{B}=\begin{bmatrix}
\til{\psi}_{\iota_1 1} & \cdots & \til{\psi}_{\iota_1 n} \\
\til{\psi}_{\iota_2 1} & \cdots & \til{\psi}_{\iota_2 n} \\
\til{\psi}_{\iota_3 1} & \cdots & \til{\psi}_{\iota_3 n} \\
\vdots & &\vdots
\end{bmatrix}
\]
we have that
\[
\sbrac{u_{ij}^G}=\til{A}\odot\til{B}.
\]
Indeed, for $s,t\iin G$ and all indices $i$ and $j$,
\[
(\til{A}\odot\til{B})_{ij}(s,t)=\inprod{\til{\xi}_i(s)}{\til{\eta}_j(t)}
=\inprod{p\xi_i(s)}{\eta_j(t)}=\inprod{\xi_i(s)}{\eta_j(t)}=u_{ij}^G(s,t).
\]
Let us identify $\bdop{\ltwo{I}}\cong\smat{I,I}$.  Then
let $\transpose{p}$ and $\wbar{p}$ denote the transpose and conjugate,
respectively, of the matrix $p$, so $\norm{\transpose{p}}=\norm{\wbar{p}}
=\norm{p}\leq 1$.
Identifying $\mat{n,I}{\linfty{G_d}}\cong\linfty{G_d,\smat{n,I}}$,
in the notation of (\ref{eq:matdecomp1}) we have that for $s\iin G$,
\[
\til{A}(s)=\begin{bmatrix} \til{\xi}_1(s) \\ \vdots \\
                           \til{\xi}_n(s) \end{bmatrix}
=\begin{bmatrix} \xi_1(s)\transpose{p} \\ \vdots \\
                 \xi_n(s)\transpose{p} \end{bmatrix}
=\begin{bmatrix} \xi_1(s) \\ \vdots \\
                 \xi_n(s) \end{bmatrix}\transpose{p}=A(s)\transpose{p}.
\]
Similarly, identifying $\mat{n,I}{\linfty{G_d}}\cong\linfty{G_d,\smat{n,I}}$,
we have that $\til{B}(s)=\wbar{p}B(s)$.  Hence $\norm{\til{A}}\leq\norm{A}$ and
$\norm{\til{B}}\leq\norm{B}$, so $\norm{\til{A}}\norm{\til{B}}\leq 1$.
Thus we see that $\sbrac{u_{ij}^G}\in\ball{\matn{\schurinv{G}}}$.  

(iv)$\implies$(i) This follows from Theorem \ref{theo:invschurmult2}.
\endpf

\smallskip
It is immediately evident from Theorem \ref{theo:haagerup} that
$\cbmg=\cbm{G_d}\cap\contb{G}$ where $G_d$ denotes the group $G$
with its discrete topology.  We observe this fact in a broader
context in Corollary \ref{cor:jsiginv2} below.

Let us consider, for the moment, the case when $G$ is amenable, and
has a left invariant mean $\m$.  We will find it notationally convenient to 
regard this mean as a finitely additive measure on $G$.

\begin{haagerup0}{\rm \cite[Prop.\ 5.2]{spronkt}}
If $G$ is amenable, then $\schurinv{G}\cong\fsalg$, completely isometrically, 
via $u\mapsto u_G$.  Moreover $\schurinv{G}$ is a completely complemented 
subspace of $\schur{G}$.
\end{haagerup0}

\proof  That $\cbmg=\fsalg$, in this case, is given in Corollary
\ref{cor:fsalgchar1}.  It remains to see that $\schurinv{G}$ is a completely 
complemented subspace of $\schur{G}$.  For $u\iin\schur{G}$ and $r\iin G$,
let $r\mult u(s,t)=u(sr,tr^{-1})$ for l.m.a.e.\ $(s,t)\iin G\cross G$.
Note that $\schurinv{G}=\{u\in\schur{G}:r\mult u=u\text{ for all }r\in G\}$.
Then define $E:\schur{G}\to\schurinv{G}$ by
\[
Eu=\weaks\int_G r\mult u d\m(r)
\]
i.e.\ $\dpair{Eu}{\mu}=\m(r\mapsto\dpair{r\mult u}{\mu})$ for all 
$\mu\iin\bloneg\htens\bloneg$.  Then it is easy to see that $E$ is
completely contractive.  Also that $E\schur{G}\subset\schurinv{G}$, and
that $Eu=u$ for $u\iin\schurinv{G}$, follow from the left invariance
of $\m$. \endpf

\smallskip
It would be interesting to determine if the (complete) complementation property
above characterises amenability of $G$.

We now return to the case of a general locally compact group $G$.
We will give a refinement of Theorem \ref{theo:haagerup}.
Let $\schurb{G}=\schurb{G,G}$, where $\schurb{G,G}$ is defined as in 
(\ref{eq:schurob}).  Let 
\begin{equation}\label{eq:schurbinv}
\schurbinv{G}=\{u\in\schurb{G}:u(sr,t)=u(s,tr^{-1})\ffor s,t,r\iin G\}.
\end{equation}
The following is a weak version of \cite{bozejkof} or of \cite{jolissaint}.

\begin{haagerup1}
\label{cor:haagerup1}
We have that $\schurinv{G}=\schurbinv{G}$.  Hence 
$\schurbinv{G}=\cbmg$ 
completely isometrically.  Thus if $u\in\cbmg$, then there is a Hilbert
space $\hil$ and $\xi,\eta\iin\contb{G,\weaks\hil}$ such that
\begin{equation}
\label{eq:schurinvip}
u(st^{-1})=\inprod{\xi(s)}{\eta(t)}
\end{equation}
for all $s,t\iin G$.  Moreover, $\vnorm{u}=\unorm{\xi}\unorm{\eta}$.
\end{haagerup1}

\proof This is immediate from the proof of
(iii)$\implies$(iv), in Theorem \ref{theo:haagerup} above 
(see (\ref{eq:tilvphitilpsi})),
and from Proposition \ref{prop:exthaagerup3}. \endpf

\smallskip
Let us finish this section by giving a second refinement of 
Theorem \ref{theo:haagerup}.  Let $\wapg$ denote the C*-algebra of 
continuous {\it weakly almost periodic} functions on $G$.  
By Grothendieck's Criterion, $\wapg$ consists of the functions
$\vphi\iin\contb{G}$ such that for any pair of sequences $\{s_n\}$ and
$\{t_m\}$ from $G$, if both of the limits $\lim_n\lim_m\vphi(s_n t_m)$ and
$\lim_m\lim_n\vphi(s_n t_m)$ exist, then they coincide.  We note that
$\wapg$ is a closed subspace of $\contb{G}$.

It is observed in \cite{xu}, that $\cbmg\subset\wapg$.  The proof is
a nice application of Grothendieck's Criterion, which we used to define
$\wapg$ above, and Corollary \ref{cor:haagerup1}.

Now we let
\[
\wschurbg=\wapg\ehtens\wapg\;\aand\;\wschurbinvg
=\wschurbg\cap\schurbinv{G}.
\]
It is clear that $\wschurbg$ is a closed subspaces of $\schurb{G}$, 
and hence $\wschurbinvg$ is a closed subspace of $\schurinv{G}$.

\begin{cbmgiswap1}
\label{cor:cbmgiswap1}
We have that $\schurinv{G}=\wschurbinvg$.  Hence $\wschurbinvg\cong\cbmg$
completely isometrically.
\end{cbmgiswap1}

\proof The functions $\til{\vphi}_{i\iota}$ and $\til{\psi}_{\iota j}$
from (\ref{eq:tilvphitilpsi}) can be approximated uniformly by sums
of translates of functions $u_{ij}$, each of which is in $\cbmg$.
Hence each $\til{\vphi}_{i\iota}$ and $\til{\psi}_{\iota j}$ is weakly
almost periodic. \endpf

\smallskip
Let us finally note that if $G$ is compact, then
\[
\schurinv{G}=\schurbinv{G}\subset\valg.
\]
where $\valg=\schurcb{G,G}=\schurco{G,G}$, 
as defined in (\ref{eq:schurcocb}).
This was shown in \cite{spronkt}.

\section{Applications}
\label{sec:applications}

\subsection{Functorial Properties of $\cbmg$}
\label{ssec:funcprop}

As a first application of the results of the previous section, we will
examine some of the functorial properties of the completely bounded
multipliers on locally compact groups.

If $X$ and $Y$ are locally compact Hausdorff spaces, and $\sig:X\to Y$ is
a continuous map, let $j_\sig:\contb{Y}\to\contb{X}$ be given by
\begin{equation}\label{eq:jsig}
j_\sig\vphi=\vphi\comp\sig.
\end{equation}
It is obvious that $j_\sig$ is a $*$-homomorphism.  Moreover, if
$\sig(X)$ is dense in $Y$ then $j_\sig$ is injective.  Indeed, if
$\vphi\in\ker j_\sig$, then $\vphi(\sig(x))=0$ for each $x\iin X$ so
$\vphi(y)=0$ for each $y\iin Y$.

Let $\schurb{X}=\schurb{X,X}$ and $\schurb{Y}=\schurb{Y,Y}$, as in
(\ref{eq:schurob}). 
Following Corollary \ref{cor:ehtensprodchar1} we get a completely contractive
algebra homomorphism $J_\sig=j_\sig\tens j_\sig:\schurb{Y}\to\schurb{X}$.
Moreover, $J_\sig$ is a complete isometry if $\sig(X)$ is dense in $Y$.

Recall that for a locally compact group $G$, $\schurbinv{G}$ is defined
in (\ref{eq:schurbinv}).

\begin{jsiginv}
\label{prop:jsiginv}
If $G$ and $H$ are locally compact groups and $\sig:G\to H$ is a
continuous homomorphism, then
$J_\sig\schurbinv{H}\subset\schurbinv{G}$ and $J_\sig$ is a complete isometry
if $\sig(G)$ is dense in $H$.
\end{jsiginv}

\proof If $u\in\schurbinv{H}$ then for $s,r,t\iin G$
\[
J_\sig u(sr,t)=u(\sig(s)\sig(r),\sig(t))=u(\sig(s),\sig(t)\sig(r^{-1}))
=J_\sig u(s,tr^{-1})
\]
so $J_\sig u\in\schurbinv{G}$.  That $J_\sig$ is a complete isometry
if $\sig(G)$ is dense in $H$ follows from comments above. \endpf

\smallskip
We obtain the main functorial property of completely bounded multipliers by
combining the above proposition with Corollary \ref{cor:haagerup1}.

\begin{jsiginv1}\label{theo:cbmfuncprop}
If $G$ and $H$ are locally compact groups and $\sig:G\to H$ is a
continuous homomorphism, then $j_\sig\cbmh\subset\cbmg$, and 
$j_\sig|_{\cbmh}:\cbmh\to\cbmg$ defines a complete contraction, which
we will again denote by $j_\sig$.  Moreover, 
if $\sig(G)$ is dense in $H$, then $j_\sig$ is a complete isometry.
\end{jsiginv1}

\smallskip
This functorial property leads to many interesting conclusions.

\begin{jsiginv2}
\label{cor:jsiginv2}
Let $G$ be a locally compact group.

{\bf (i)} The inclusion $\cbmg\hookrightarrow\cbm{G_d}$ is a complete
isometry, where $G_d$ denotes the group $G$ with discrete topology.

{\bf (ii)} If $\galp$ denotes the almost periodic compactification of $G$ and
$\sigalp:G\to\galp$ the canonical map {\rm (see \cite[Sec.\ 6.1]{dixmierB})}, 
then $j_{\sigalp}:\fal{\galp}\to\cbmg$ is a complete isometry.

{\bf (iii)} If $H$ is a closed subgroup of $G$, then the restriction map
$u\mapsto u|_H$ is a complete contraction from $\cbmg$ to $\cbmh$.
Moreover, if $H$ is open, then this map is a surjective complete quotient map.
If $H$ is amenable, then this map is a complete quotient map onto its range.

{\bf (iv)}  If $N$ is a closed normal subgroup of $G$, and $q:G\to G/N$ is the 
quotient map, then $j_q:\cbm{G/N}\to\cbmg$ is a complete isometry.
\end{jsiginv2}

\proof  For (i), the inclusion is just the map $j_{\iota_d}$, 
where $\iota_d:G_d\to G$
is the identity map.  For (ii), we just note that $\galp$ is compact,
so $\fal{\galp}=\fsal{\galp}=\cbm{\galp}$ completely isometrically, by
Corollary \ref{cor:fsalgchar1}.  

For (iii), the restriction map is just
$j_\iota$, where $\iota:H\to G$ is the inclusion map.  
If $H$ is open, the result  
follows from the equivalence of (i) and (ii) in Theorem \ref{theo:haagerup}.
If $H$ is amenable,
then $j_\iota|_{\fsalg}:\fsalg\to\fsalh$ is a complete quotient map onto
its range.
Indeed, the range of $j_\iota|_{\fsalg}$ is $\ccfs{\vpi|_H}$ where
$\vpi$ is the universal representation of $G$.  The adjoint of
$j_\iota|_{\fsalg}$ is then the inclusion
$\vn{\vpi|_H}\hookrightarrow\wstarg$, and hence is a complete isometry.  
We thus obtain the factorization 
$j_\iota|_{\fsalg}=j_\iota\comp i$, where $i:\fsalg\to\cbmg$ is
the completely contractive inclusion map.  Hence $j_\iota$ must be a complete
quotient onto its range.

For (iv), $j_q$ is complete isometry, since $q$ is a surjection. \endpf

We remark that in (iv), above, the range of $j_q$ consists of elements
of $\cbmg$ which are constant on cosets of $N$.  In \cite[p.\ 682]{haagerupk}
(and following that, \cite[Sec.\ 2.2]{xu}), it is claimed that
every element of $\cbmg$ which is constant on cosets of $N$ is in
the range of $j_q$.  If
$G$ is amenable, this result follows from that $\cbmg=\fsalg$
and \cite[2.26]{eymard}. 

Regarding (iii), we note that it is known for an
amenable group $G$ that it is possible to have a closed subgroup $H$
for which the restriction map from $\fsalg$ to $\fsalh$ is not surjective
(see \cite[p.\ 204]{eymard}, for example), though it is always
a complete quotient onto its range.  We are unable to determine
if, in general, the restriction map from $\cbmg$ to $\cbmh$ is
a (complete) quotient map onto its range.

It follows from (ii) above that if $\pi$ is a finite dimensional unitary
representation of $G$, then the imbedding $\ccfs{\pi}\hookrightarrow
\cbmg$ is a complete isometry.  This result is generalised below.

\begin{jsiginv3}
\label{cor:jsiginv3}
Let $G$ be a locally compact group and $\pi$ be a continuous unitary
representation such that there is an amenable group $H$, a continuous
homomorphism $\sig:G\to H$ and a continuous unitary representation 
$\pi'$ of $H$, such that $\pi=\pi'\comp\sig$.  Then 
the imbedding $\ccfs{\pi}\hookrightarrow\cbmg$ is a
complete isometry.
\end{jsiginv3}

\proof First note that $j_\sig:\fsal{H}\to\fsalg$ is a complete contraction,
since its adjoint $j_\sig^*:\wstarg\to\wstarh$ is a $*$-homomorphism,
by the universal property of $\wstarg$.  Note too that each $u\iin\ccfs{\pi}$ 
is of the form $j_\sig u'$ for some $u'\iin\ccfs{{\pi'}}\subset\fsal{H}$
by \cite[Prop.\ 2.10]{arsac}, and
recall that $\fsal{H}=\cbmh$ completely isometrically, since $H$ is amenable.
Without loss of generality we may assume that $\sig(G)$ is dense in $H$, so  
$j_\sig:\cbmh\to\cbmg$ is a complete isometry.
Then if $\sbrac{u_{ij}}\in\matn{\ccfs{\pi}}$, we have that
$\sbrac{u_{ij}}=\nlift{j_\sig}\sbrac{u_{ij}'}$ for some
$\sbrac{u_{ij}'}\iin\matn{\fsal{H}}$, and hence
\[
\cbmnorm{\sbrac{u_{ij}}}\leq\norm{\sbrac{u_{ij}}}
\leq\norm{\sbrac{u'_{ij}}}
=\cbmnorm{\sbrac{u'_{ij}}}
=\cbmnorm{\nlift{j_\sig}\sbrac{u'_{ij}}}
=\cbmnorm{\sbrac{u_{ij}}}.
\]
Thus $\cbmnorm{\sbrac{u_{ij}}}=\norm{\sbrac{u_{ij}}}$. \endpf

\subsection{A Predual for $\cbmg$}

In this section, we give a concrete construction of a predual of 
$\cbmg$.  The existence of this predual was recognised in \cite{decanniereh},
and in a more general form in \cite{krausr}, where it was identified as an
operator predual.

Note that $\sig:\schur{G}\to\schur{G}$ given by
$\sig u(s,t)=u(s,t^{-1})$ is a bijective complete isometry.
Thus we may now consider
\begin{equation}
\sig\schurinv{G}=\cbrac{u\in\schur{G}:
\begin{matrix} u(sr,t)=u(s,rt)\text{ for all }r\iin  \\
G\text{ and l.m.a.e. }(s,t)\in G\cross G\end{matrix}}.
\end{equation}
As in (\ref{eq:usubg}), given $u\iin\sig\schurinv{G}$, we can define
a function $u_G:G\to\Cee$ such that for l.a.e.\ $t\iin G$
\[
u_G(t)=u(s,s^{-1}t)
\]
for l.a.e.\ $s\iin G$.
Recalling (\ref{eq:lonehdual}), we have that $\schur{G}\cong
(\blone{G}\htens\blone{G})^*$, where the dual pairing is given by
\[
\dpair{u}{f\tens g}=\int_G \int_G u(s,t)f(s)g(t)dsdt
\]
for each elementary tensor $f\tens g\iin\bloneg\htens\bloneg$.
If $\fX\subset\bloneg\htens\bloneg$, let $\fX^\perp=\{u\in\schur{G}:
\dpair{u}{x}=0\text{ for all }x\iin\fX\}$.

Let $\m_0:\bloneg\phtens\bloneg\to\bloneg$ be the multiplication map.
This is an unbounded map on an incomplete space.  Let $\fK$ be the closure
of $\ker\m_0$ in $\bloneg\htens\bloneg$ and
\begin{equation}
\queg=\brac{\bloneg\htens\bloneg}/\fK.
\end{equation}
Endow $\queg$ with the quotient operator space structure.

\begin{cbmpredual}\label{theo:cbmpredual}
$\sig\schurinv{G}=\fK^\perp$, and hence $\queg^*\cong\schurinv{G}\cong\cbmg$ 
as operator spaces.
\end{cbmpredual}

We will call $\queg$ {\it the standard predual} of $\cbmg$.  Although
its uniqueness is not assured, we will refer to the
topology $\sig(\cbmg,\queg)$ as {\it the} weak* topology on $\cbmg$.

\smallskip
\proof To see that $\sig\schurinv{G}\subset\fK^\perp$, it suffices to see that
$\sig\schurinv{G}\subset(\ker\m_0)^\perp$.  To that end, if 
$u\in\sig\schurinv{G}$
and $\mu=\sum_{k=1}^n f_k\tens g_k\iin\ker\m_0$, then
\begin{align*}
\dpair{u}{\mu}&=\sum_{k=1}^n\int_G \int_G u(s,t)f_k(s)g_k(t) dt ds \\
&=\sum_{k=1}^n\int_G \int_G u(s,s^{-1}t)f_k(s)g_k(s^{-1}t) dt ds \\
&=\sum_{k=1}^n\int_G u_G(t)f_k\con g_k(t)dt=\int_G u_G(t)\m_0(\mu)(t)=0.
\end{align*}
Conversely, if $u\in(\ker\m_0)^\perp$, $f, g\in\bloneg$ and $\del_r$ is
the Dirac measure at $r$ for $r\iin G$, then 
$\m_0(f\con\del_r\tens g-f\tens\del_r\con g)=0$, so
\begin{align*}
0&=\dpair{u}{f\con\del_r\tens g-f\tens\del_r\con g} \\
&=\int_G \int_G u(s,t)
  \brac{ \frac{1}{\Del(r)}f(sr^{-1})g(t)-f(s)g(r^{-1}t) }dtds \\
&=\int_G \int_G (u(sr,t)-u(s,rt))f(s)g(t)dtds
\end{align*}
and hence $u(sr,t)=u(s,rt)$ for marginally almost every $(s,t)\iin G\cross G$,
so $u\in\sig\schurinv{G}$.  Thus $(\ker\m_0)^\perp\subset\sig\schurinv{G}$.
\endpf

\smallskip
The map $\m_0:\bloneg\phtens\bloneg\to\bloneg$ is surjective by Cohen's
Factorization Theorem.  For $f\iin\bloneg$ let
\[
\qnorm{f}=\inf\{\hnorm{\mu}:\mu\in\bloneg\phtens\bloneg\aand
\m_0(\mu)=f\}.
\]
The function $\qnorm{\cdot}$ is clearly a semi-norm on $\bloneg$.
Since $\bloneg$ has a contractive bounded approximate identity,
the Cohen-Hewitt Factorization Theorem \cite[(32.22)]{hewittrII} 
gives for any $\eps>0$ that
we can write $f=g\con h$, where $\norm{g}_1\leq 1$ and $\norm{f-h}_1<\eps$,
so that $\norm{g}_1\norm{h}_1<\norm{f}_1+\eps$.  Hence $\qnorm{f}\leq
\norm{f}_1$.  Let $\m:\bloneg\htens\bloneg\to\queg$ be the quotient map.
If $f=g\con h$ in $\bloneg$, let $\qmap(f)=\m(g\tens h)$.  Then
$\qmap:\bloneg\to\queg$ is clearly well-defined and 
linear with $\norm{\qmap(f)}=\qnorm{f}$.

The following shows that we may consider $\queg$ to be the completion
of $\bloneg$ under the norm $\qnorm{\cdot}$.  Hence our space
$\queg$ is isomorphic to the predual discovered in \cite{decanniereh}.

\begin{cbmpredual1}
The map $\qmap:\bloneg\to\queg$ is an injective (completely) contractive
map with dense range.  Moreover, for $f\iin\bloneg$ and $u$ in $\cbmg$, we have
$\dpair{u}{\qmap(f)}=\int_G u(s)f(s)ds$.
\end{cbmpredual1}

\proof That $\qmap$ is contractive was shown above.  That $\qmap(\bloneg)$
is dense in $\queg$
follows from that $\bloneg\phtens\bloneg$ is dense in $\bloneg\htens\bloneg$.
If $u\in\cbmg$, let $u^G(s,t)=u(st)$ so $u\in\sig\schurinv{G}$. 
Then for $f\iin\bloneg$, writing $f=g\con h$, we have that
\begin{align*}
\dpair{u}{\qmap(f)}&=\dpair{u^G}{g\tens h}
=\int_G\int_G u(st)g(s)h(t)dtds \\
&=\int_G\int_G u(t)g(s)h(s^{-1}t)dtds
=\int_G u(t)g\con h(t)dt
=\int_G u(t)f(t)dt.
\end{align*}
That $\qmap$ is injective follows from that $\falg$ is both contained in 
$\cbmg$ and dense in $\conto{G}$, so
for every $f\iin\bloneg$ there is $u\iin\falg$ such that $\int_G f(t)u(t)dt
\not=0$. \endpf

\smallskip
We now want to see how a continuous unitary representation of $G$ induces a map
on $\queg$.  There is an alternative approach in \cite{kraus}.
The author is grateful to J.\ Kraus for showing him this paper.

If $\pi$ is a continuous unitary representation of $G$, let $\pi_1$
be its extension to $\bloneg$.  We have that $\pi_1\tens\pi_1:
\bloneg\htens\bloneg\to\cstaro{\pi}\htens\cstaro{\pi}$ is a
complete contraction.  Let $\m_h:\cstaro{\pi}\htens\cstaro{\pi}\to\cstaro{\pi}$
be the completely contractive multiplication map. Then let
\begin{equation}\label{eq:multpi}
\m_{\pi}=\m_h\comp(\pi_1\tens\pi_1):\bloneg\htens\bloneg\to\cstaro{\pi}.
\end{equation}
Since $\m_{\pi}$ is continuous and
$\pi_1\comp\m_0=\m_\pi$ on $\bloneg\phtens\bloneg$,
$\ker\m_\pi\supset\wbar{\ker\m_0}$.  Hence $\m_\pi$ induces
a map $\vqmap{\pi}:\queg=\bloneg\htens\bloneg/\wbar{\ker\m_0}\to\cstaro{\pi}$
by
\begin{equation}\label{eq:vqmap}
\vqmap{\pi}(\m(x))=\m_\pi(x)
\end{equation}
for $x\iin\bloneg\htens\bloneg$, where $\m:\bloneg\htens\bloneg\to\queg$ is 
the quotient map.  Note that $\vqmap{\pi}^*:\wcfs{\pi}\to\cbmg$ is the
completely contractive inclusion map, so $\vqmap{\pi}$ is completely 
contractive.  
Finally, note that the diagram
\begin{equation}\label{eq:multpiandvamap}
\xymatrix{
\bloneg\htens\bloneg \ar[dr]_{\m} \ar[rr]^<<<<<<<<<<<<<<{\m_\pi} & &
\cstaro{\pi} \\
& \queg \ar[ur]_{\vqmap{\pi}} 
}
\end{equation}
commutes.

It is natural to ask when $\vqmap{\pi}$ is injective, since it is
tempting to want to think of $\queg$ as a dense subspace of $\cstarg$.
However, there is a technical obstruction against obtaining this.
Recall that $\vpi$ denotes the universal representation of $G$,
and $\fsalg=\ccfs{\vpi}$.

\begin{vmapinjective}
\label{prop:vmapinjective}
The map
$\vqmap{\pi}$ is injective if and only if $\ccfs{\pi}$ is weak* dense
in $\cbmg$.  In particular, $\vqmap{\vpi}:\queg\to\cstarg$ is injective
if and only if $\fsalg$ is weak* dense in $\cbmg$.
\end{vmapinjective}

\proof $\vqmap{\pi}$ is injective if and only if 
$(\cstaro{\pi})^*\cong\wcfs{\pi}$ is point separating for $\queg$.
By the Hahn-Banach Theorem, this happens exactly when $\wcfs{\pi}$
is weak* dense in $\cbmg$.  Since $\ccfs{\pi}$ is 
$\sig(\wcfs{\pi},\cstaro{\pi})$-dense in
$\wcfs{\pi}$, $\ccfs{\pi}$ is $\sig(\wcfs{\pi},\vqmap{\pi}(\queg))$-dense
in $\wcfs{\pi}$, and thus
its weak* closure coincides with that of $\wcfs{\pi}$ in $\cbmg$. \endpf

\smallskip
We say that $G$ has the {\it approximation property} 
\cite{haagerupk} if $\falg$ is weak* dense in $\cbmg$.  
Recall that $\lam$ denotes the left regular representation of $G$, and 
$\falg=\ccfs{\lam}$.

\begin{vmapinjective1}
$\vqmap{\lam}$ is injective if and only if $G$ satisfies the
approximation property.  
\end{vmapinjective1}

It is natural to ask when the maps $\m_\pi$ and $\vqmap{\pi}$ are surjective.

\begin{qmapsurj}\label{prop:qmapsurj}
For a continuous unitary representation $\pi$ of a
locally compact group $G$, the following are equivalent.

{\bf (i)} The map $\m_\pi:\bloneg\htens\bloneg\to\cstaro{\pi}$ is surjective.

{\bf (ii)} The map $\vqmap{\pi}:\queg\to\cstaro{\pi}$ is surjective.

{\bf (iii)} The injection $\wcfs{\pi}\hookrightarrow\cbmg$
is bounded below.
\end{qmapsurj}

\proof That (i) is equivalent to (ii) follows from (\ref{eq:multpiandvamap}).
That (ii) is equivalent to (iii) follows from Theorem
\ref{theo:cbmpredual} and the Open Mapping Principle.  \endpf

\smallskip
In what follows
we  obtain an interesting characterisation of the amenability of $G$
which is essentially a reformulation of the fact that
\begin{equation}\label{eq:losert}
\begin{matrix} \text{\it if the injection }\falg\hookrightarrow\cbmg \\
\text{\it is bounded below, then }G\;\text{\it is amenable.}\end{matrix}
\end{equation}
This highly non-trivial fact was proved by V.\ Losert, likely in
a similar vein as \cite{losert84}, though it remains
unpublished.  The author is grateful to Z.-J.\ Ruan for showing him this.  
It follows then that $\fsalg=\cbmg$ only if $G$ is amenable.
This fact is known for discrete groups by the attractive paper 
\cite{bozeko85}.  
In \cite{losert84} it was shown that $\mults\falg=\fsalg$ implies
that $G$ is amenable.  If $G$ has the approximation property, then
a nice proof that $\fsalg=\cbmg$ implies that $G$ is amenable
can be found in \cite{krausr99}. 

The characterisation below
can also be found in \cite{pisier98} in the case that $G$ is discrete.

\begin{gamenchar}
The following are equivalent for any locally compact group $G$:

{\bf (i)}  $G$ is amenable.

{\bf (ii)} The map $\m_\vpi:\bloneg\htens\bloneg\to\cstarg$ is surjective.

{\bf (iii)} The map $\m_\lam:\bloneg\htens\bloneg\to\cstarrg$ is surjective.
\end{gamenchar}

\proof (i)$\implies$(ii)  If $G$ is amenable, then $\fsalg=\cbmg$ 
isometrically.  In particular, condition (iii) of Proposition 
\ref{prop:qmapsurj} is met.

The implication (ii)$\implies$(iii) follows from the commutativity of the 
following diagram.
\[
\xymatrix{
\bloneg\htens\bloneg \ar[r]^{\vpi_1\tens\vpi_1}
\ar[dr]_{\lam_1\tens\lam_1}
&\cstarg\htens\cstarg \ar[r]^>>>>>>{\m_h} 
&\cstarg \ar[d]^{\lam_*} \\
&\cstarrg\htens\cstarrg \ar[r]^>>>>>>{\m_h} & \cstarrg
}
\]

That (iii)$\implies$(i) follows from Proposition \ref{prop:qmapsurj} and
(\ref{eq:losert}). \endpf
 
\smallskip
We note that since $\bloneg$ has a bounded approximate identity, the condition
(ii) above is equivalent to saying that $\bloneg$ has {\it similarity degree 
2}, in the terminology of \cite{pisier98}.  
Furthermore, this condition implies the following
result, obtained in the discrete case by Pisier in \cite{pisier98}, and 
announced by him for the case of non-discrete groups in \cite{pisierP}:

\begin{simprob}\label{theo:simprob}
A locally compact group $G$ is amenable if and only if for 
every bounded representation $\pi:G\to\bdophinv$ (i.e.\ for which 
$\unorm{\pi}=\sup\{\norm{\pi(s)}:s\in G\}<+\infty$), there is
an operator $S\iin\bdophinv$ such that $S\pi(s)S^{-1}$ is unitary
for each $s\iin G$, and $\norm{S}\norm{S^{-1}}\leq\unorm{\pi}^2$.
\end{simprob}

\noindent The necessity condition above is from \cite{dixmier}.
Details of Theorem \ref{theo:simprob} are verified in \cite{spronkTh}.

{
\bibliography{biblio}
\bibliographystyle{plain}
}

\smallskip
{\sc Department of Mathematics, Texas A\& M University, College Station, 
Texas 77843-3368, U.S.A.}

{\it E-mail address:} {\tt spronk@math.tamu.edu}

\end{document}